\documentclass[10pt,a4paper]{article}
\usepackage{amsmath}
\usepackage{amssymb}
\usepackage{amsthm}
\usepackage{amsfonts}
\usepackage{amscd}
\usepackage[mathscr]{eucal}
\usepackage{multirow}
\newcommand{\Z} {{\mathbb  Z}}
\newcommand{\Q}{{\mathbb  Q}}
\newcommand{\F}{{\mathbb  F}}
\newcommand{\C}{{\mathbb  C}}

\begin{document}
\parindent  25pt
\baselineskip  10mm
\textwidth  15cm    \textheight  23cm \evensidemargin -0.06cm
\oddsidemargin -0.01cm

\title{ {On Twists of A Family of Elliptic Curves and Their $ L-$Function
 }}
\author{\mbox{}
{ Derong Qiu }
\thanks{ \quad E-mail:
derong@mail.cnu.edu.cn, \ derongqiu@gmail.com } \\
(School of Mathematical Sciences,
 Capital Normal University, \\
Beijing 100048, P.R.China )  }

\date{}
\maketitle
\parindent  24pt
\baselineskip  10mm
\parskip  0pt

\par   \vskip 0.4cm

{\bf Abstract} \quad Let $ E $ be an elliptic curve defined over a number field,
the conjecture of Birch and Swinnerton-Dyer (BSD, for short) asserts
a deep relation between the group $ E(K) $ of rational points
and the $ L-$function $ L(E/K, s)$ of $ E $ at $ s = 1. $
Very few explicit results about $ E(K) $ and $ L(1) $ are known, even no general method
is known to determine $ L(1) $ vanishing or not for a given elliptic curve.
In this paper, we study some quantities related to BSD of a special class of elliptic curves,
more precisely, we study the arithmetic of quadratic twists
of elliptic curves $ y^{2} = x(x + \varepsilon p )(x +
  \varepsilon q) $ and their $L-$function. Based on some classical works, especially those
of Greenberg, Kramer-Tunnell, Kato-Rohrlich, Manin and Mazur, under some conditions, we obtain results
about the vanishing of the value at $ s = 1 $ of the $ L$-function, and explicitly
determine the following quantities: the norm index $ \delta (E, \Q, K), $ the root numbers,
the set of anomalous prime numbers, a few prime numbers at which the image of Galois
representation are surjective. We also study the relation between the ranks of the
Mordell-Weil groups, Selmer groups and Shafarevich-Tate groups, and
the structure about the $ l^{\infty }-$Selmer groups and the Mordell-Weil groups
over $ \Z_{l}-$extension via Iwasawa theory. These results provide some useful
evidence toward verifying the BSD for a family of elliptic curves.
\par  \vskip  0.2 cm

\par  \vskip  0.1 cm


\par     \vskip  0.4 cm

\hspace{-0.8cm}{\bf 1. Introduction }

\par \vskip 0.2 cm

Let $ E $ be an elliptic curve over a number field $ K, $ and $ L(E/K, s) $ be the
$ L-$function of $ E $ over $ K. $ By Mordell-Weil theorem
(see, e.g. [Sil1]), the set $ E(K) $ of $ K-$rational points of $ E $ is a finitely
generated abelian group. Hence
$$ E(K) \simeq \Z^{r} \bigoplus E(K)_{\text{tors}}, $$
where $ r = \text{rank}(E(K)) \geq 0 $ is the rank of $ E $ over $ K, $ and
$ E(K)_{\text{tors}} $ is the torsion subgroup of $  E(K). $
\par \vskip 0.2 cm

{\bf Conjecture 1.1} (see [Sil1]). \ The $ L-$function $ L(E/K, s) $ of $ E $ over $ K $ has an analytic
continuation to the entire complex plane, and satisfies a functional equation relating
the values at $ s $ and $ 2-s. $   \\
This conjecture was proved when $ K = \Q $ (see [BCDT], [TW], [Wi]).  \\
The conjecture of Birch and Swinnerton-Dyer (BSD, for short)
for elliptic curves states that
\par \vskip 0.2 cm

{\bf Conjecture 1.2} (Birch and Swinnerton-Dyer conjecture, see [Sil1]). \\
(1) \ The rank of $ E(K) $ equals the order of vanishing of $ L(E/K, s) $ at $ s = 1. $  \\
(2)  $$ \lim_{s\rightarrow 1} \frac{L(E/K, s)}{(s-1)^{r}} =
\frac{\Omega \cdot \sharp \amalg\hskip-7pt\amalg(E/K) \cdot R(E/K)\cdot
\prod_{v| \mathcal{N}} c_{v}(E)}{\sharp E(K)_{\text{tors}}^{2}}, $$
where $ r = \text{rank}E(K), \Omega = $ the real period, $ E(K)_{\text{tors}} $ is the torsion
subgroup of $ E(K), R(E/K) $ is the regulator of $ E(K)/E(K)_{\text{tors}},  \mathcal{N} $ is
the conductor of $ E/K,  c_{v}(E) = [E(K_{v}) : E_{0}(K_{v})] $ is the Tamagawa number of $ E $
at the place $ v, \ \amalg\hskip-7pt\amalg(E/K) $ is the Shafarevich-Tate group of $ E $
over $ K, $ which is conjectured to be a finite group.
\par \vskip 0.2 cm

In the literature, much important progress has been made about the BSD conjecture. For example, 
for elliptic curves over the rational number field $ \Q, $ let $ r_{an}(E/\Q) $ denote the order 
of vanishing of $ L(E/\Q, s) $ at $ s = 1. $ Then one current state of the BSD conjecture is expressed
by the result:
\par \vskip 0.2 cm

{\bf Theorem 1.3} (Gross-Zagier, Kolyvagin, etc., see [Kol3]). The equality
$ \text{rank}E(\Q) = r_{an}(E/\Q) $ holds and
$ \sharp \amalg\hskip-7pt\amalg(E/\Q) $ is finite if $ r_{an}(E/\Q) \leq 1.$  \\
Yet, at present, to explicitly determine the arithmetic quantities such as $ E(K) $ and the order
of $ L(E/K, s) $ at $ s = 1 $ are generally not easy, even for the question about determining
whether the value $ L(E/K, 1) $ vanishing or not.   \\
In this paper, we will study explicitly $ L(1) $ and some related arithmetic quantities about
twists of a family of elliptic curves $ E $ over the rational number field $ \Q, $ from which,
for example, we obtain that $ L(E_{d} / \Q, 1) = 0 $ for many quadratic twists $ E_{d} $ of $ E. $
More precisely, we consider the elliptic curves $$
  E = E^{\varepsilon }: \ y^{2} = x(x + \varepsilon p )(x +
  \varepsilon q), \quad (\varepsilon = \pm 1), \eqno(1.1) $$ and
their quadratic $ D-$twist
$$  E_{D}= E^{\varepsilon }_{D}:  \ y^{2} = x(x + \varepsilon p D)(x +
  \varepsilon q D), \quad \eqno(1.2) $$
where $ p $ and $ q $ are odd prime numbers with $ q - p = 2, $ and
$ D = D_{1} \cdots D_{n} $ is a square-free integer with distinct
odd prime numbers $ D_{1}, \cdots, D_{n} $ satisfying $ (pq, D) = 1.
$ When $ D = 1,  E_{1} = E, $ and for $ \varepsilon = 1 $ (resp. $ -1$),
we sometimes write $ E^{\varepsilon } = E^{+} $ (resp. $ E^{-}$). By Tate's
algorithm (see [Ta], [Sil2]), the discriminant, $j-$invariant and
conductor of $ E_{D}/ \Q $ are obtained as follows, respectively
$$ \Delta = 64p^{2}q^{2}D^{6}, \ j = \frac{64(p^{2} + 2q)^{3}}
{p^{2}q^{2}}, \ N_{E_{D}} = 2^{5}pqD^{2}.  \eqno(1.3) $$
So the equation (1.2) above is a global minimal Weierstrass equation
for $ E_{D} $ over the rational number field $ \Q. $ Moreover,
$ E_{D} / \Q $ has additive reduction at $ 2, D_{1}, \cdots, D_{n}, $
has multiplicative reduction at $ p, q, $ and has good reduction
at other finite places. \\
In the following, we study the arithmetic of these
elliptic curves. The following quantities are explicitly
determined: the norm index $ \delta (E, \Q, K)$ (see Theorem 3.3),
the root numbers (see Theorem 5.3), the set of anomalous prime numbers
(see Proposition 2.4), a few prime numbers at which the
image of Galois representation are surjective (see Proposition 2.7).
The relation between the ranks of the Mordell-Weil
groups, Selmer groups and Shafarevich-Tate groups, and the structure
about the $ l^{\infty }-$Selmer groups and the Mordell-Weil groups
over $ \Z_{l}-$extension via Iwasawa theory are studied (see
Propositions 3.1, 4.1, 4.2, and Theorems 3.4, 3.7, 3.8, 4.3, 4.4).
On $ L(1), $ one of our main result is as follows
\par \vskip 0.2 cm

{\bf Theorem 1.4} (see Theorem 5.5 below) \ Let $ E = E^{\varepsilon } $ be the
elliptic curve in (1.1) and let $ K = \Q(\sqrt{\mu D}) $ be the quadratic number
field with $ D $ in (1.2) and $ \mu = \pm 1. $ We assume that
$ D \equiv \mu \ (\text{mod} \ 4). $ Let $ L(E / \Q, s) =
\Sigma _{n = 1}^{\infty } a_{1}(n) n^{-s} $ be the $ L-$function
as above. Let $ E_{\mu D} / \Q $  be the quadratic $ (\mu D)-$twist
of $ E / \Q, $ and $ \chi _{K} $ be the quadratic Dirichlet character
associated to $ K. $ \\
(1) \ Assume one of the following two hypotheses holds: \\
(a) \ $ \varepsilon = 1 $ and $ p \equiv 5, 7 \ (\text{mod} \ 8); $ \\
(b) \ $ \varepsilon = -1 $ and $ p \equiv 3, 5 \ (\text{mod} \ 8). $ \\
Then $ L(E / \Q, 1) = 2 \Sigma _{n = 1}^{\infty } \frac{a_{1}(n)}{n}
e^{-n\pi /2\sqrt{2pq}}. $ \\
further, for all integer $ r \geq 0, $
\begin{align*}
&L^{(r)}(E / \Q, 1) = 2 \pi \Sigma _{n = 1}^{\infty } a_{1}(n)
\int _{1 /4 \sqrt{2pq}}^{\infty } [\log ^{r} t +
(-1)^{r}\log ^{r}(2^{5}pq t)] e^{-2n\pi t} dt. \ \text{also}, \\
&L(E_{\mu D} / \Q, 1) = (1 + \chi _{K}(-2pq)) \cdot
\Sigma _{n = 1}^{\infty } \frac{a_{1}(n)}{n} \chi _{K}(n) \cdot
e^{-n\pi /2D\sqrt{2pq}},
\end{align*}
In particular, if $ \chi _{K}(-2pq) = -1, $ then
$ L(E_{\mu D} / \Q, 1) = 0. $ \\
(2) \ Assume one of the following two hypotheses holds: \\
(a$^{\prime}$) \ $ \varepsilon = 1 $ and $ p \equiv 1, 3 \ (\text{mod} \ 8); $ \\
(b$^{\prime}$) \ $ \varepsilon = -1 $ and $ p \equiv 1, 7 \ (\text{mod} \ 8). $ \\
Then $ L(E / \Q, 1) = 0, $ \\
further, for all integer $ r \geq 0, $
\begin{align*}
&L^{(r)}(E / \Q, 1) = 2 \pi \Sigma _{n = 1}^{\infty } a_{1}(n)
\int _{1 /4 \sqrt{2pq}}^{\infty } [\log ^{r} t +
(-1)^{r + 1}\log ^{r}(2^{5}pq t)] e^{-2n\pi t} dt. \ \text{also}, \\
&L(E_{\mu D} / \Q, 1) = (1 - \chi _{K}(-2pq)) \cdot
\Sigma _{n = 1}^{\infty }  \frac{a_{1}(n)}{n} \chi _{K}(n) \cdot
e^{-n\pi /2D\sqrt{2pq}}.
\end{align*}
In particular, if $ \chi _{K}(-2pq) = 1, $ then $ L(E_{\mu D} / \Q, 1) = 0. $  \\
(For some concrete example on $ L(1), $ see Example 5.6 below).
\par \vskip 0.2 cm

These results, together with some former results about Mordell-Weil
groups and Selmer groups as in [QZ1] and [LQ], provide some useful
evidence toward verifying the BSD for a family of elliptic curves, which we will
discuss in a separate paper.
\par     \vskip  0.2 cm

{\it Organisation of the paper.} \ Section 2 includes some basic facts on reduction
from Tate's algorithm, and some results on anomalous prime, ramification and Galois
representation deduced from the works of Mazur, Bahargava-Skinner-Zhang and Serre.
In Section 3, by using Kramer's method and Kramer-Tunnell' formula, and the former results
in [Q1], [QZ1], we compute the norm index, Tamagawa number, Selmer group, rank, and 
some congruences between rank and Shafarevich-Tate group. In Section 4, following 
mainly the works of Mazur, Greenberg and Kato-Rohrlich, we study the structure about 
the $ l^{\infty }-$Selmer groups and the Mordell-Weil groups over $ \Z_{l}-$extension via 
Iwasawa theory. Finally, in Section 5, by results of Rohrlich, we compute the root numbers, 
and by using a formula of Manin on $ L(1), $ we obtain some results on the vanishing of 
the value at $ s = 1 $ of the $ L$-function.

\par     \vskip  0.4 cm

\hspace{-0.8cm}{\bf 2. Reduction, ramification and Galois
representation}

\par \vskip 0.2 cm

In the following, unless otherwise stated, every conclusion for
the elliptic curves $ E_{D} $ in (1.2) also holds for
$ E_{1} = E $ in (1.1) when take $ D = 1. $ For a prime number
$ l $ and an integer $ m, \ (\frac{m}{l}) $ is the usual Legendre
quadratic residue symbol.
\par \vskip 0.2 cm

{\bf Lemma 2.1} \ Let $ E_{D} / \Q $ be the elliptic curve in (1.2)
above. \\
(1) \ At each prime $ l \mid N_{E_{D}}, $ the Kodaira type is as follows: \\
$ III $ for $ l = 2; \ I_{2} $ for $ l = p $ or $ q; $ and $ I_{0}^{\ast }
$ for $ l = D_{1}, \cdots, D_{n}, $ respectively. \\
The Tamagawa number $ c_{l} $ is as follows: \\
$ c_{l} = 2 $ for $ l = 2, p, q; $ and $ c_{l} = 4 $ for
$ l = D_{1}, \cdots, D_{n}. $ \\
(2) \ $ E_{D} $ has split multiplicative reduction at $ p $ if and only
if $ (\frac{2\varepsilon D}{p}) = 1. $ \\
(3) \ $ E_{D} $ has split multiplicative reduction at $ q $ if and only
if $ (\frac{- 2 \varepsilon D}{q}) = 1. $ \\
(4) \ Let $ l $ be a prime number such that $ l \nmid 2pqD. $ Then
$ E_{D} $ has good supersingular reduction at $ l $ if and only
if $ \sum _{m = 0}^{(l-1)/2} \begin{pmatrix} \frac{l - 1}{2} \\
m \end{pmatrix}^{2} p^{m}
q^{\frac{l-1}{2} - m} \equiv 0 \ (\text{mod} \ l). $ \\
(5) \ The torsion subgroup $ E_{D}(\Q)_{\text{tors}} \cong
\Z / 2 \Z \times \Z / 2 \Z, $ and for $ D = 1, $ we have
$ E(F)_{\text{tors}} \cong \Z / 2 \Z \times \Z / 2 \Z $ for any
quadratic number field $ F. $ \\
(6) \ Assume $ 3 \nmid pqD. $ Let $ F $ be a number field, and let
$ \mathfrak{p} $ be a prime ideal of $ F $ lying over $ 3, $ let
$ e = e(\mathfrak{p} / 3) $ and $ f = f(\mathfrak{p} / 3) $ be the
ramification index and residue degree, respectively. Then we have \\
(6a) \ if $ e(\mathfrak{p} / 3) = f(\mathfrak{p} / 3) = 1, $ then
$ E_{D}(F)_{\text{tors}} \cong \Z / 2 \Z \times \Z / 2 \Z; $ \\
(6b) \ if $ f(\mathfrak{p} / 3) = 1 $ and $ E_{D} $ has additive
reduction at some finite places of $ F $ lying over $ 2, $ then
$ E_{D}(F)_{\text{tors}} \cong \Z / 2 \Z \times \Z / 2 \Z $ or
$ \Z / 2 \Z \times \Z / 6 \Z; $ \\
(6c) \ if $ f(\mathfrak{p} / 3) = 1, $ then
$ E_{D}(F)_{\text{tors}} / E_{D}(F)[3^{\infty }] \cong
\Z / 2 \Z \times \Z / 2 \Z, $ where $ E_{D}(F)[3^{\infty }] $
denotes the $ 3-$primary component of $ E_{D}(F)_{\text{tors}}; $ \\
(6d) \ If $ E_{D} $ has an additive reduction at some finite places
of $ F $ lying over $ 2, $ then $ \sharp E_{D}(F)_{\text{tors}} = 2^{m} $
or $ 2^{m} \cdot 3 $ for some $ m \in \Z_{\geq 0}. $
\par  \vskip 0.1 cm

{\bf Proof.} \ (1) is a consequence of direct calculation by
the Algorithm of [Ta]; (2), (3) and (4) are easily obtained
(see [Sil1] for the methods); (5) follows from Lemma 2 and Lemma 4
of [QZ2]; (6) is similar to the Prop.1 in [QZ1, p.1374].
\quad $ \Box $ \\
Particularly, by (2) and (3) of Lemma 2.1, one can easily see that,
$ E^{+} $ has split multiplicative reduction at both $ p $ and $ q $ if
$ p \equiv 1, 7 \ (\text{mod} \ 8), $ and has non-split multiplicative reduction
at both $ p $ and $ q $ if $ p \equiv 3, 5 \ (\text{mod} \ 8); $ Also,
$ E^{-} $ has split multiplicative reduction at $ p $ and non-split multiplicative reduction
at $ q $ if $ p \equiv 1, 3 \ (\text{mod} \ 8), $ and has non-split
multiplicative reduction at $ p $ and split multiplicative reduction at $ q $ if
$ p \equiv 5, 7 \ (\text{mod} \ 8). $
\par  \vskip 0.2 cm

{\bf Corollary 2.2.} \ For the elliptic curves $ E_{D} / \Q $ in (1.2)
above, \\
(1) \ $ E_{D} $ has good supersingular reduction at $ 3 $ if
$ 3 \nmid pqD; $ \\
(2) \ $ E_{D} $ has good ordinary reduction at $ 5 $ if
$ 5 \nmid pqD; $ \\
(3) \ $ E_{D} $ has good ordinary reduction at $ 7 $ if
$ 7 \nmid pqD $ and $ p \equiv 1, 4 \ (\text{mod} \ 7); $ \\
(4) \ $ E_{D} $ has good supersingular reduction at $ 7 $ if
$ 7 \nmid pqD $ and $ p \equiv 2, 3, 6 \ (\text{mod} \ 7). $
\par  \vskip 0.1 cm

{\bf Proof.} \ Follows easily from the above Lemma 2.1(4).
\quad $ \Box $

\par  \vskip 0.2 cm
For an elliptic curve $ E / \Q $ and a prime number $ l, $ we
denote the reduction of $ E $ at $ l $ by $ \widetilde{E}_{l}, $
and let $ a_{l} = l + 1 - \sharp \widetilde{E}_{l}(\F_{l}), $
where $ \F_{l} $ is the field with $ l $ elements. For a positive
integer $ m, \ E[m] = \{P \in E(\overline{\Q}) : \ m P = 0 \} $
is the group of $ m-$division points of $ E, $ where
$ \overline{\Q} $ is an algebraic closure of $ \Q. $ Let
$ G_{\Q} = \text{Gal}(\overline{\Q} / \Q) $ be the
absolute Galois group, and let
$ \rho _{l} : \ G_{\Q} \longrightarrow \
\text{Gl}_{2}(\F_{l}) $ be the Galois representation of
$ G_{\Q} $ given by the action of $ G_{\Q} $ on the $ l-$division
points of $ E $ (see, e.g., [Sil1, p.90]). By the open
image theorem of Serre ([Se1]), $ \rho _{l}$ is surjective
for all but finitely many prime numbers $ l. $
\par  \vskip 0.2 cm

{\bf Lemma 2.3.} \ For the elliptic curves $ E_{D} / \Q $ in (1.2)
above, \\
(1) \ if $ 3 \nmid pqD, $ then
$ \sharp \widetilde{E_{D,}}_{3}(\F_{3}) = 4 $ and $ a_{3} = 0. $ \\
(2) \ if $ 7 \nmid pqD, $ and $ p \equiv 2, 3, 6 \ (\text{mod} \ 7), $
then $ \sharp \widetilde{E_{D,}}_{7}(\F_{7}) = 8 $ and
$ a_{7} = 0. $ \\
(3) \ assume $ 5 \nmid pqD, $  \\
(3a) \ if $ p \equiv 1, 2 \ (\text{mod} \ 5), $ then
$$ \sharp \widetilde{E_{D,}}_{5}(\F_{5}) =
\left \{\begin{array}{l} 4 \quad \text{if}
\ D \equiv 1, 4 \ (\text{mod} \ 5) \\
8 \quad \text{if} \ D \equiv 2, 3 \ (\text{mod} \ 5),
\end{array} \right. \ \text{and} \ a_{5} =
\left \{\begin{array}{l} 2 \quad \text{if}
\ D \equiv 1, 4 \ (\text{mod} \ 5) \\
-2 \quad \text{if} \ D \equiv 2, 3 \ (\text{mod} \ 5),
\end{array} \right. $$
(3b) \ if $ p \equiv 4 \ (\text{mod} \ 5), $ then
$$ \sharp \widetilde{E_{D,}}_{5}(\F_{5}) =
\left \{\begin{array}{l} 8 \quad \text{if}
\ D \equiv 1, 4 \ (\text{mod} \ 5) \\
4 \quad \text{if} \ D \equiv 2, 3 \ (\text{mod} \ 5),
\end{array} \right. \ \text{and} \ a_{5} =
\left \{\begin{array}{l} -2 \quad \text{if}
\ D \equiv 1, 4 \ (\text{mod} \ 5) \\
2 \quad \text{if} \ D \equiv 2, 3 \ (\text{mod} \ 5).
\end{array} \right. $$
(4) \ assume $ 7 \nmid pqD, $  \\
(4a) \ if $ \left \{\begin{array}{l} \varepsilon = 1 \\
p \equiv 1 \ (\text{mod} \ 7)
\end{array} \right. \ \text{or} \
\left \{\begin{array}{l} \varepsilon = -1 \\
p \equiv 4 \ (\text{mod} \ 7),
\end{array} \right. $
then
$$ \sharp \widetilde{E_{D,}}_{7}(\F_{7}) =
\left \{\begin{array}{l} 12 \quad \text{if}
\ D \equiv 1, 2, 4 \ (\text{mod} \ 7) \\
4 \quad \text{if} \ D \equiv 3, 5, 6 \ (\text{mod} \ 7),
\end{array} \right. \ \text{and} \ a_{7} =
\left \{\begin{array}{l} -4 \quad \text{if}
\ D \equiv 1, 2, 4 \ (\text{mod} \ 7) \\
4 \quad \text{if} \ D \equiv 3, 5, 6 \ (\text{mod} \ 7),
\end{array} \right. $$
(4b) \ if $ \left \{\begin{array}{l} \varepsilon = 1 \\
p \equiv 4 \ (\text{mod} \ 7)
\end{array} \right. \ \text{or} \
\left \{\begin{array}{l} \varepsilon = -1 \\
p \equiv 1 \ (\text{mod} \ 7),
\end{array} \right. $
then
$$ \sharp \widetilde{E_{D,}}_{7}(\F_{7}) =
\left \{\begin{array}{l} 4 \quad \text{if}
\ D \equiv 1, 2, 4 \ (\text{mod} \ 7) \\
12 \quad \text{if} \ D \equiv 3, 5, 6 \ (\text{mod} \ 7),
\end{array} \right. \ \text{and} \ a_{7} =
\left \{\begin{array}{l} 4 \quad \text{if}
\ D \equiv 1, 2, 4 \ (\text{mod} \ 7) \\
-4 \quad \text{if} \ D \equiv 3, 5, 6 \ (\text{mod} \ 7).
\end{array} \right. $$
(5) \ $ \sharp \widetilde{E_{D,}}_{2}(\F_{2}) = 3, \quad
\sharp \widetilde{E_{D,}}_{_{D_{i}}}(\F_{D_{i}}) = D_{i} + 1 \
(i = 1, \cdots, n), $
$$ \sharp \widetilde{E_{D,}}_{p}(\F_{p}) =
\left \{\begin{array}{l} p \quad \quad \ \text{if}
\ (\frac{2 \varepsilon D}{p}) = 1 \\
p + 2 \quad \text{if} \ (\frac{2 \varepsilon D}{p}) = -1,
\end{array} \right. \ \text{and} \
\sharp \widetilde{E_{D,}}_{q}(\F_{q}) =
\left \{\begin{array}{l} q \quad \quad \ \text{if}
\ (\frac{-2 \varepsilon D}{q}) = 1 \\
q + 2 \quad \text{if} \ (\frac{-2 \varepsilon D}{q}) = -1.
\end{array} \right. $$
{\bf Proof.} \ Via direct calculation.
\quad $ \Box $
\par  \vskip 0.2 cm

Recall that a prime number $ l $ is said to be anomalous for
an elliptic curve $ E / \Q $ if $ E $ has good reduction at $ l $
and $ \sharp \widetilde{E}_{l}(\F_{l}) \equiv 0 \ (\text{mod} \ l) $
(see [Ma2, p.186] and [M, p.25]). We denote $ \text{Anom}(E / \Q) =
\{l : \ l \ \text{is an anomalous prime number for} \ E / \Q \}. $
\par  \vskip 0.2 cm

{\bf Proposition 2.4.} \ For the elliptic curves $ E_{D} / \Q $ in (1.2)
above, we have $ \text{Anom}(E_{D} / \Q) = \emptyset . $
\par  \vskip 0.1 cm

{\bf Proof.} \ Since the conductor $ N_{E_{D}} = 2^{5}pqD^{2}, $ we have
$ \ 2, p, q, D_{i} \notin \text{Anom}(E_{D} / \Q) \ (i = 1, \cdots, n). $
On the other hand, by Lemma 2.1(5) above, $ E_{D}(\Q)_{\text{tors}} \cong
\Z / 2 \Z \times \Z / 2 \Z, $ so by the results 2.10(b) of [M, p.26] we
have $ \text{Anom}(E_{D} / \Q) \subset \{2, 3, 5 \}, $ and so
$ \text{Anom}(E_{D} / \Q) \subset \{3, 5\}. $ For $ l = 3 $ or $ 5, $
we may assume that $ l \nmid pqD, $ then by Lemma 2.3(1) and (3) above,
we have $ \sharp \widetilde{E_{D,}}_{3}(\F_{3}) = 4 $ and
$ \sharp \widetilde{E_{D,}}_{5}(\F_{5}) = 4 $ or $ 8, $ which shows that
$ 3, 5 \notin \text{Anom}(E_{D} / \Q), $ so
$ \text{Anom}(E_{D} / \Q) = \emptyset . $ \quad $ \Box $
\par  \vskip 0.2 cm

For our next discussion, we need the following
\par  \vskip 0.2 cm

{\bf Lemma 2.5} (see [BSZ, p.4] and [Sil2, Prop.6.1 and exer.V.5.13]).
\ Let $ E $ be an elliptic curve over $ \Q $ with conductor $ N_{E}. $
Let $ l, l^{\prime } $ be two prime numbers with $ l \neq l^{\prime }. $
Suppose $ l \parallel N_{E}. $ Then $ E[l^{\prime }] $ is ramified at
$ l $ if and only if $ l^{\prime } \nmid \ \text{ord}_{l}(\Delta _{l}) $
for a minimal discriminant $ \Delta _{l} $ of $ E $ at $ l. $
\par  \vskip 0.2 cm

{\bf Proposition 2.6.} \ For the elliptic curves $ E_{D} / \Q $ in
(1.2) above , let $ l $ be a prime number. Then \\
(1) \ $ E_{D}[l] $ is ramified at $ p $ if and only if $ l > 2 $
and $ l \neq p; $ \\
(2) \ $ E_{D}[l] $ is ramified at $ q $ if and only if $ l > 2 $
and $ l \neq q. $ \\
In particular, $ E_{D}[p] $ is ramified at $ q, $ and $ E_{D}[q] $ is
ramified at $ p. $
\par  \vskip 0.1 cm

{\bf Proof.} \ Since the equation in (1.2) above is global minimal
for $ E_{D} / \Q, $ we have
$ \Delta _{l} = \Delta = 64p^{2}q^{2}D^{6} $ for any prime number
$ l, $ so $$ \text{ord}_{l}(\Delta _{l}) =
\left \{\begin{array}{l} 0 \quad \text{if}
\ l \nmid 2pqD  \\
6 \quad \text{if} \ l \mid 2D  \\
2 \quad \text{if} \ l = p \ \text{or} \ q.
\end{array} \right. $$
On the other hand, the conductor $ N_{E_{D}} = 2^{5}pqD^{2}, $ so a prime
number $ l \parallel N_{E_{D}} \Leftrightarrow l = p $ or $ q. $ By the above
discussion, $ \text{ord}_{p}(\Delta _{p}) =
\text{ord}_{q}(\Delta _{q}) = 2, $ so the conclusion follow from the
above Lemma 2.5. \quad $ \Box $
\par  \vskip 0.2 cm

{\bf Proposition 2.7.} \ For the elliptic curves $ E_{D} / \Q $ in
(1.2) above, let $ l $ be a prime number, and $ \rho _{l} $ be the
corresponding Galois representation. \\
(1) \ If $ 3 \nmid pqD, $ then $ \rho _{3} $ is surjective,
i.e., $ \rho _{3}(G_{\Q}) = \text{Gl}_{2}(\F_{3}). $  \\
(2) \ If $ 7 \nmid pqD $ and $ p \equiv 2, 3, 6 \ (\text{mod} \ 7), $
then $ \rho _{7} $ is surjective, i.e., $ \rho _{7}(G_{\Q}) =
\text{Gl}_{2}(\F_{7}). $  \\
(3) \ If $ 3 \nmid pqD, \ l \nmid pqD $ and $ l > 3105, $ then
$ \rho _{l} $ is surjective, i.e., $ \rho _{l}(G_{\Q})
= \text{Gl}_{2}(\F_{l}). $
\par  \vskip 0.1 cm

{\bf Proof.} \ (1) \ Under the assumption, by Cor.2.2(1) above,
$ E_{D} $ has good supersingular reduction at $ 3; $ also, the
discriminant $ \Delta = (2D)^{6}(pq)^{2} $ is obviously not a cube,
so the conclusion follows from Serre's theorem (see [Se1] or
[PR, Prop.4.4]). \\
(2) \ Under the assumption, by Cor.2.2(4) above,
$ E_{D} $ has good supersingular reduction at $ 7; $ also, since
the conductor $ N_{E_{D}} = 2^{5}pqD^{2} $ and the invariant
$ j = \frac{64(p^{2} + 2q)^{3}}{p^{2}q^{2}}, $ we have
$ p \parallel N $ and $ \text{ord}_{p}(j) = -2 \not\equiv 0 \
(\text{mod} \ 7). $ So the conclusion follows from Serre's theorem
(see [Se1] or [PR, Prop.4.4]). \\
(2) \ Under the assumption, $ 3 $ is the smallest (odd) prime number
at which $ E_{D} $ has good reduction. Also, $ j \notin \Z $ and
$ \text{ord}_{p}(j) = -2 < 0. $ Moreover, the prime number $ l $ under
our assumption obviously satisfies $ l > (\sqrt{3} + 1)^{8}. $ So
the conclusion follows from Prop.24 of [Se1]. \quad $ \Box $

\par     \vskip  0.4 cm

\hspace{-0.8cm}{\bf 3. Rank, norm index, Shafarevich-Tate group
and $ l-$Selmer group}

\par  \vskip 0.2 cm

Let $ E / \Q $ be the elliptic curve in (1.1) above, and let $ K =
\Q (\sqrt{D}) $ be the quadratic number field, where
$ D = D_{1} \cdots D_{n} $ with distinct odd prime numbers
$ D_{1}, \cdots, D_{n} $ as in (1.2) above. Let
$ M_{K} $ be a complete set of places on $ K, $ and $ M_{K}^{\infty } $
(resp. $ M_{K}^{0} $) its subset of infinite (resp. finite) places.
 Let $ S_{K} = M_{K}^{\infty } \cup \{v \in M_{K}^{0} : \ v \mid 2pq \}. $
The group of $ S_{K}-$units of $ K $ is denoted by $ U_{K,S}, $ the
ideal class group of $ K $ is denoted by $ \text{Cl}(K), $ and the
$ S_{K}-$class group of $ K $ is denoted by $ \text{Cl}_{S}(K), $
precisely, $ \text{Cl}_{S}(K) $ is the quotient of $ \text{Cl}(K) $
by the subgroup generated by the classes represented
by the finite primes in $ S_{K} $ (see [Sa, p.127]). For an abelian
group $ A $ and a positive integer $ m, $ we write $ A[m] =
\{a \in A : \ m a = 0 \}. $ For a vector space $ V $ over $ \F_{2}, $
we denote its dimension by $ \text{dim}_{2} V. $ For a finitely
generated abelian group $ A, $ we denote its rank by
$ \text{rank}(A). $ The next result is about $ E(K), $ the group of
rational points of $ E $ over $ K. $
\par  \vskip 0.2 cm

{\bf Proposition 3.1.} \ Let $ E / \Q $ be the elliptic curve in
(1.1), and $ K = \Q (\sqrt{D}) $ be the quadratic number field
as above, we have $ \text{rank}(E(K)) \leq 14 +
2 \text{dim}_{2} \text{Cl}_{S}(K)[2]. $
\par  \vskip 0.1 cm

{\bf Proof.} \ Let $ E^{\prime } : \ y^{2} = x^{3} -
2\varepsilon (p + q) x^{2} + 4x. $ There is an
isogeny $ \varphi $ of degree $ 2 $ between $ E $ and $ E^{\prime } $
with the dual isogeny $ \widehat{\varphi } $ as in
[QZ1, pp.1372,1373]. Let $ \text{Sel}_{\varphi }(E/K) $ and
$ \text{Sel}_{\widehat{\varphi } }(E^{\prime }/K) $ be the
$ \varphi -$Selmer group of $ E/K $ and the
$ \widehat{\varphi }-$Selmer group of $ E^{\prime }/K, $
respectively, and $ \amalg\hskip-7pt\amalg(E/K) $ (resp. $
\amalg\hskip-7pt\amalg(E^{\prime }/K) $ be the Shafarevich-Tate
groups of $ E/K $ (resp. $ E^{\prime }/K $) (see [Sil1, Chapt.10]).
Then (see [Sil1, pp298, 301])
\begin{align*} &\text{dim}_{2} E(K)/ 2E(K) + \text{dim}_{2}
E^{\prime }(K)[\widehat{\varphi }] / \varphi (E(K)[2]) \\
&= \text{dim}_{2} \text{Sel}_{\varphi }(E/K) -
\text{dim}_{2} \amalg\hskip-7pt\amalg(E/K) [\varphi] +
\text{dim}_{2} \text{Sel}_{\widehat{\varphi } }(E^{\prime }/K)
- \text{dim}_{2} \amalg\hskip-7pt\amalg(E^{\prime }/K)
[\widehat{\varphi }]. \end{align*}
Note that $ E^{\prime }(K)[\widehat{\varphi }] = \{O, (0, 0) \}, \
\varphi (E(K)[2]) = \{O, (0, 0) \}, $ so $ \text{rank}( E(K)) \leq
\text{dim}_{2} \text{Sel}_{\varphi }(E/K) +
\text{dim}_{2} \text{Sel}_{\widehat{\varphi } }(E^{\prime }/K) - 2. $
On the other hand, the following exact sequence
is known (see, e.g., [St, p.5], [Sz, p.55]): \
$ 0 \rightarrow U_{K,S} / U_{K,S}^{2} \rightarrow K(S_{K}, 2)
\rightarrow \text{Cl}_{S}(K)[2] \rightarrow 0, $ where,
$ K(S_{K}, 2) = \{b K^{\ast ^{2}} \in K^{\ast } / K^{\ast ^{2}}: \
\text{ord}_{v}(b) \equiv 0 \ (\text{mod} \ 2) \
\text{for all} \ v \notin S_{K} \}. $ So by the Dirichlet unit
theorem (see [L, pp.104, 105]), we have $ \text{dim}_{2} K(S_{K}, 2)
= \sharp S_{K} + \text{dim}_{2} \text{Cl}_{S}(K)[2] \leq 8 +
\text{dim}_{2} \text{Cl}_{S}(K)[2] $ because $ \sharp S_{K} =
\sharp M_{K}^{\infty } + \sharp \{v \in M_{K}^{0} : \ v \mid 2pq \}
\leq 2 + 6 = 8. $ Also, $ \sharp \text{Sel}_{\varphi }(E/K)
\leq \sharp K(S_{K}, 2) $ and
$ \sharp \text{Sel}_{\widehat{\varphi } }(E^{\prime }/K)
\leq \sharp K(S_{K}, 2) $ (see [Sil1, p.302]), so from the above discussion,
$ \text{rank}( E(K)) \leq 2 \text{dim}_{2} K(S_{K}, 2) - 2 \leq
14 + 2 \text{dim}_{2} \text{Cl}_{S}(K)[2]. $ \quad $ \Box $
\par  \vskip 0.2 cm

Next, we need state some notations. Let $ F $ be a number field and $ L $
be a quadratic extension of $ F, $ we write $ M_{F} $ (resp.$
M_{L} $) for a complete set of places on $ F $ (resp.$ L $). Fix a
place $ w \in M_{L} $ lying above $ v $ for each $ v \in M_{F}. $
Denote the Galois group $ \text{Gal}(L_{w} / F_{v}) $ by $ G_{w}, $
where $ F_{v} $ and $ L_{w} $ are the completions of $ F $ at $ v $
and $ L $ at $ w, $ respectively. Let $ E $ be an elliptic curve
over $ F. $ For every $ v \in M_{F}, $ we denote
$ \delta _{v} = \log _{2}(E(F_{v}) : N(E(L_{w}))), $ this is the local
norm index studied deeply in [Kr] and [KT]. For some of their arithmetic
application (see,e.g., [MR], [Q1]). Let $ \delta (E, F, L) $ be the
sum of all the local norm index, i.e., $ \delta (E, F, L) =
\sum _{v \in M_{F}} \delta _{v}. $ Now, for the elliptic
curve $ E / \Q $ in (1.1) and the quadratic number field
$ K = \Q(\sqrt{D}) $ as above, we come to calculate explicitly
the quantity $ \delta (E, \Q, K) $ as in [Q1, p.5054, and Section 3
there], and give some application.
\par  \vskip 0.2 cm

{\bf Lemma 3.2.} \ Let $ E / \Q $ be the elliptic curve in
(1.1), $ \mu = \pm 1, $ and $ K = \Q (\sqrt{\mu D}) $ be the
quadratic number field with square-free integer
$ D = D_{1} \cdots D_{n} $ as in (1.2)
above. Fix a place $ w \in M_{K} $ lying above $ 2. $ Let
$ \Delta _{w}, c_{w} $ and $ f_{w} $ be the minimal discriminant,
Tamagawa number and the exponent of the conductor of $ E $ at $ w $
(i.e., over $ K_{w}$)(see [Sil1]), respectively. \\
(1) \ If $ D \equiv 5 \mu (\text{mod} 8), $ then
$ K_{w} \cong \Q_{2}(\sqrt{-3}), $ and \\
Type $ III, \ \text{ord}_{w}(\Delta _{w}) = 6, \ f_{w} = 5, $ and
$ c_{w} = 2. $ \\
(2) \ If $ D \equiv 7 \mu (\text{mod} 8), $ then
$ K_{w} \cong \Q_{2}(\sqrt{-1}), $ and \\
Type $ I_{2}^{\ast }, \ \text{ord}_{w}(\Delta _{w}) = 12, \ f_{w} = 6, $
and $ c_{w} = \left \{\begin{array}{l} 2 \quad \text{if}
\ p \equiv 1 (\text{mod} 4)  \\
4 \quad \text{if} \ p \equiv 3 (\text{mod} 4).
\end{array} \right. $  \\
(3) \ If $ D \equiv 3 \mu (\text{mod} 8), $ then
$ K_{w} \cong \Q_{2}(\sqrt{3}), $ and \\
Type $ I_{2}^{\ast }, \ \text{ord}_{w}(\Delta _{w}) = 12, \ f_{w} = 6, $
and $ c_{w} = \left \{\begin{array}{l} 4 \quad \text{if}
\ p \equiv 1 (\text{mod} 4)  \\
2 \quad \text{if} \ p \equiv 3 (\text{mod} 4).
\end{array} \right. $
\par  \vskip 0.1 cm

{\bf Proof.} \ For the case $ \mu D \equiv 3, 5, 7 (\text{mod} 8), $
from the proof of Lemma3.1 in [Q1, p.5057], we have \
$ K_{w} \cong \Q_{2}(\sqrt{-3}) \Longleftrightarrow \mu D \equiv
5 \ (\text{mod} \ 8); \ K_{w} \cong \Q_{2}(\sqrt{-1})
\Longleftrightarrow \mu D \equiv 7 \ (\text{mod} \ 8); \
K_{w} \cong \Q_{2}(\sqrt{3}) \Longleftrightarrow \mu D \equiv
3 \ (\text{mod} \ 8). $ Then the conclusion follows from Tate's
algorithm (see [Ta], [Sil2]), in a way as done in the proof of
Lemma3.1 of [Q1, p.5057]. \quad $ \Box $
\par  \vskip 0.2 cm

{\bf Theorem 3.3.} \ Let $ E / \Q $ be the elliptic curve in
(1.1), $ \mu = \pm 1, $ and $ K = \Q (\sqrt{\mu D}) $ be the
quadratic number field with square-free integer
$ D = D_{1} \cdots D_{n} $ as in (1.2)
above. Denote $ \mu _{0} = (1 - \mu ) / 2. $ Then we have
$ 2 n + \mu _{0} \leq \delta (E, \Q, K) \leq 2 n + 4 + \mu _{0}. $
More precisely,  \\
(1) \ $ \delta (E, \Q, K) = 2 n + \mu _{0} $ if and only if
$ D \equiv \mu (\text{mod} 8) $ and $ (\frac{\mu D}{p}) =
(\frac{\mu D}{q}) = 1. $ \\
(2) \ $ \delta (E, \Q, K) = 2 n + 1 + \mu _{0} $ if and only if one of the
following four hypotheses holds : \\
(2a) \ $ D \equiv 5 \mu (\text{mod} 8) $ and
$ (\frac{\mu D}{p}) = (\frac{\mu D}{q}) = 1; $  \\
(2b) \ $ D \equiv 7 \mu (\text{mod} 8), \ p \equiv 3 (\text{mod} 4) $
and $ (\frac{\mu D}{p}) = (\frac{\mu D}{q}) = 1;  $ \\
(2c) \ $ D \equiv 3 \mu (\text{mod} 8), \ p \equiv 1 (\text{mod} 4) $
and $ (\frac{\mu D}{p}) = (\frac{\mu D}{q}) = 1; $ \\
(2d) \ $ D \equiv \mu (\text{mod} 8) $ and $ (\frac{\mu D}{p}) +
(\frac{\mu D}{q}) = 0. $  \\
(3) \ $ \delta (E, \Q, K) = 2 n + 2 + \mu _{0} $ if and only if one of the
following six hypotheses holds: \\
(3a) \ $ D \equiv 7 \mu (\text{mod} 8), \ p \equiv 1 (\text{mod} 4) $
and $ (\frac{\mu D}{p}) = (\frac{\mu D}{q}) = 1; $ \\
(3b) \ $ D \equiv 3 \mu (\text{mod} 8), \ p \equiv 3 (\text{mod} 4) $
and $ (\frac{\mu D}{p}) = (\frac{\mu D}{q}) = 1; $ \\
(3c) \ $ D \equiv 5 \mu (\text{mod} 8) $ and $ (\frac{\mu D}{p}) +
(\frac{\mu D}{q}) = 0; $ \\
(3d) \ $ D \equiv 7 \mu (\text{mod} 8), \ p \equiv 3 (\text{mod} 4) $
and $ (\frac{\mu D}{p}) + (\frac{\mu D}{q}) = 0; $ \\
(3e) \ $ D \equiv 3 \mu (\text{mod} 8), \ p \equiv 1 (\text{mod} 4) $
and $ (\frac{\mu D}{p}) + (\frac{\mu D}{q}) = 0; $ \\
(3f) \ $ D \equiv \mu (\text{mod} 8) $ and $(\frac{\mu D}{p}) =
(\frac{\mu D}{q}) = -1. $  \\
(4) \ $ \delta (E, \Q, K) = 2 n + 3 + \mu _{0} $ if and only if
one of the following five hypotheses holds: \\
(4a) \ $ D \equiv 7 \mu (\text{mod} 8), \ p \equiv 1 (\text{mod} 4) $
 and $ (\frac{\mu D}{p}) + (\frac{\mu D}{q}) = 0; $ \\
(4b) \ $ D \equiv 3 \mu (\text{mod} 8), \ p \equiv 3 (\text{mod} 4) $
and $ (\frac{\mu D}{p}) + (\frac{\mu D}{q}) = 0; $  \\
(4c) \ $ D \equiv 5 \mu (\text{mod} 8) $ and
$ (\frac{\mu D}{p}) = (\frac{\mu D}{q}) = -1; $  \\
(4d) \ $ D \equiv 7 \mu (\text{mod} 8), \ p \equiv 3 (\text{mod} 4) $
and $ (\frac{\mu D}{p}) = (\frac{\mu D}{q}) = -1; $ \\
(4e) \ $ D \equiv 3 \mu (\text{mod} 8), \ p \equiv 1 (\text{mod} 4) $
and $ (\frac{\mu D}{p}) = (\frac{\mu D}{q}) = -1. $ \\
(5) \ $ \delta (E, \Q, K) = 2 n + 4 + \mu _{0} $ if and only
if one of the following two hypotheses holds: \\
(5a) \ $ D \equiv 7 \mu (\text{mod} 8), \ p \equiv 1 (\text{mod} 4) $
and $ (\frac{\mu D}{p}) = (\frac{\mu D}{q}) = -1; $  \\
(5b) $ D \equiv 3 \mu (\text{mod} 8), \ p \equiv 3 (\text{mod} 4) $
and $ (\frac{\mu D}{p}) = (\frac{\mu D}{q}) = -1. $
\par  \vskip 0.1 cm

{\bf Proof.} \ We consider the case $ \mu = 1, $ the other case is similar.
Let $ S $ be the set of finite places of $ \Q $
obtained by collecting together all places that ramify in $ K / \Q $
and all places of bad reduction for $ E / \Q, $ so
$ S = \{2, p, q, D_{1} \cdots D_{n} \}. $ Although the cases here become
more complicated, we will take our calculation in a way as in the
Lemma 3.2 of [Q1, p.5058], so we need to use the same notations $ S_{0},
S_{g}, S_{gu}, S_{ar}, S_{a}, S_{smr}, S_{nsmr},  S_{nsmr}^{\prime },
S_{nsmr}^{\prime \prime} $ as in the Remark of [Q1, pp.5055,5056].
For the convenience of the reader, we write them in the present case as: \\
$ S_{0} = \{v \in S: \ v \ \text{is ramified or inertial in} \ K \}; \\
S_{g} = \{v \in S_{0}: \ v \nmid 2 \ \text{and} \ E \
\text{has good reduction at} \ v \} = \{D_{1}, \cdots, D_{n} \}; \\
S_{gu} = \{v \in S_{0}: v \mid 2, \ E \ \text{has good reduction
at} \ v \ \text{and} \ \Q_{v} \ \text{is unramified over} \
\Q_{2} \} \\
= \emptyset ; \\
S_{ar} = \{v \in S_{0}: \ E \ \text{has additive reduction at} \ v
\} = \left \{\begin{array}{l} \{2 \} \quad \text{if}
\ D \equiv 3, 5, 7 (\text{mod} 8)  \\
\emptyset \quad \text{if} \ D \equiv 1 (\text{mod} 8);
\end{array} \right.  \\
S_{a} = S_{ar} \cup \{v \in S_{0}: v \mid 2, \ E \ \text{has good
reduction at} \ v \ \text{and} \ \Q_{v} \
\text{is ramified over} \ \Q_{2} \} \\
= S_{ar}; \\
S_{smr} = \{v \in S_{0}: \ E \ \text{has split multiplicative
reduction at} \ v \} \subset \{p, q \} \cap S_{0}; \\
S_{nsmr} = \{v \in S_{0}: \ E \ \text{has non-split
multiplicative reduction at} \ v \} \\
= S_{nsmr}^{\prime } \sqcup S_{nsmr}^{\prime \prime} \ (\text{the
disjoint union})  \subset \{p, q \} \cap S_{0}, \quad
\text{where} \\
S_{nsmr}^{\prime } = \{v \in S_{nsmr}: \ v \ \text{is inertial
in} \ K \} = S_{nsmr}, \\
S_{nsmr}^{\prime \prime} = \{v \in S_{nsmr}: \ v \ \text{is
ramified in} \ K \} = \emptyset . $ \\
Obviously, $ S_{0} = S_{g} \sqcup S_{gu} \sqcup S_{a} \sqcup S_{smr}
\sqcup S_{nsmr} $ (the disjoint union). \\
By definition, $ \delta (E, \Q, K) =
\sum _{v \in M_{\Q}} \delta _{v}, $ where $ \delta _{v} =
\log _{2}(E(\Q_{v}) : N(E(K_{w}))) $ is the local norm index.
Furthermore, by the results in [Kr], one can obtain that
$ \delta (E, \Q, K) =
\delta _{\infty } + \delta _{f}, $ where $ \delta _{\infty } $ is
as in [Q, p.5054], and $ \delta _{f} =
\delta _{g} + \delta _{m} + \delta _{a} $ with
$ \delta _{g}, \delta _{m}, \delta _{a} $ in [Q1, pp.5055,5056],
that is,
\begin{align*}
&\delta _{a} = \sum _{v \in S_{a}} \delta _{v}; \quad
\delta _{m} = \delta _{smr} + \delta _{nsmr} \ \text{with} \ \delta
_{smr} = \frac{1}{2} \sum _{v \in S_{smr}} \left(1 + (\Delta _{v},
D)_{\Q_{v}} \right) \ \text{and} \\
&\delta _{nsmr} = \frac{1}{2} \sum _{v \in S_{nsmr}^{\prime }}
\left(1 + (-1)^{v (\Delta _{v})} \right) + \sum _{v \in
S_{nsmr}^{\prime \prime}} \left( \frac{1}{2} \left(1 + (\Delta _{v},
D)_{\Q_{v}} \right) \cdot (-1)^{v (\Delta _{v})} + 1 \right); \\
&\delta _{g} = \sum _{v \in S_{g}} \dim _{2}
\widetilde{E_{v}}(k_{v}) [2] + \sum _{v \in S_{gu}} \varepsilon (v), \quad
\text{where}
\end{align*}
 $$ \varepsilon (v) = \left \{\begin{array}{l}
\frac{1}{2} \left(1 - (-1)^{v(D)} \right) \cdot [\Q_{v} : \Q_{2}] \
\text{if} \ E \ \text{has good supersingular reduction at} \ v, \\
\\
\frac{1}{2} (3 + (\Delta _{v}, D)_{\Q_{v}}) \quad \quad \text{if} \ E
\ \text{has good ordinary reduction at} \ v.
\end{array}
\right. $$ Here $ \widetilde{E}_{v} $ is the reduction of $ E $ at $
v, \ k_{v} $ is the residue field of $ \Q_{v}, $ and
$ ( , )_{\Q_{v}} $ is the Hilbert symbol
(see [Se 2, Chapt.XIV]). \\
It is easy to see here that $ \delta _{\infty } = 0 $
since $ D > 0. $ So we only need to calculate
$ \delta _{g}, \delta _{m},  \delta _{a}. $ For this, we divide our
discussion into the following cases. \\
Case for $ \delta _{g}. $ \ Since $ E $ has good reduction at each
$ D_{i} (i = 1, \cdots, n), $
we have an injective homomorphism $ E(\Q)_{\text{tors}} \hookrightarrow
\widetilde{E}_{D_{i}}(\F_{D_{i}}) $ (see [Kn, p.130]). So by
Lemma 2.1(5) above, we have $ \widetilde{E}_{D_{i}}(\F_{D_{i}})[2]
\cong (\Z/ 2 \Z)^{2}. $ and so \\
$ \delta _{g} = \sum _{l \in S_{g}} \dim_{2}
\widetilde{E}_{l}(\F_{l})[2] = \sum _{i = 1}^{n} \dim_{2}
\widetilde{E}_{D_{i}}(\F_{D_{i}})[2] = 2 n, $
i.e., $ \delta _{g} = 2 n. $ \\
Case for $ \delta _{m}. $ \ Since the equation (1.1) is global minimal
for $ E / \Q, $ we have $ \text{ord}_{p}(\Delta _{p}) =
\text{ord}_{q}(\Delta _{q}) = 2, $ so
$ 1 + (-1)^{\text{ord}_{l}(\Delta _{l})} = 2 $ for $ l = p $ or $ q, $
and so $ \delta _{nsmr} = \sharp S_{nsmr}. $ Also $
(\Delta _{p}, D)_{\Q_{p}} = (\Delta _{q}, D)_{\Q_{q}} = 1 $ because
$ \Delta _{p} = \Delta _{q} = (8pq)^{2}. $ So
$ \delta _{smr} = \sharp S_{smr}. $ Hence $ \delta _{m} =
\sharp S_{smr} + \sharp S_{nsmr} = \sharp (S_{0} \cap \{p, q \})
\leq 2. $ The set $ S_{0} $ can be determined as follows. \\
If $ D \equiv 1 (\text{mod} 8), $ then $ S_{0} =
\left \{\begin{array}{l} \{D_{1}, \cdots, D_{n}, p \} \quad \text{if}
\ (\frac{D}{p}) = -1 \ \text{and} \ (\frac{D}{q}) = 1  \\
\{D_{1}, \cdots, D_{n}, q \} \quad \text{if}
\ (\frac{D}{p}) = 1 \ \text{and} \ (\frac{D}{q}) = -1   \\
\{D_{1}, \cdots, D_{n} \} \quad \text{if}
\ (\frac{D}{p}) = (\frac{D}{q}) = 1   \\
\{D_{1}, \cdots, D_{n}, p, q \} \quad \text{if}
\ (\frac{D}{p}) = (\frac{D}{q}) = -1;
\end{array} \right.  $  \\
If $ D \equiv 3, 5, 7 (\text{mod} 8), $ then $ S_{0} =
\left \{\begin{array}{l} \{2, D_{1}, \cdots, D_{n}, p \} \quad
\text{if} \ (\frac{D}{p}) = -1 \ \text{and} \ (\frac{D}{q}) = 1  \\
\{2, D_{1}, \cdots, D_{n}, q \} \quad \text{if}
\ (\frac{D}{p}) = 1 \ \text{and} \ (\frac{D}{q}) = -1   \\
\{2, D_{1}, \cdots, D_{n} \} \quad \text{if}
\ (\frac{D}{p}) = (\frac{D}{q}) = 1   \\
\{2, D_{1}, \cdots, D_{n}, p, q \} \quad \text{if}
\ (\frac{D}{p}) = (\frac{D}{q}) = -1.
\end{array} \right.  $  \\
From this, we get \\
$ \delta _{m} =
\left \{\begin{array}{l} 0 \quad \text{if}
\ (\frac{D}{p}) = (\frac{D}{q}) = 1   \\
1 \quad \text{if} \ (\frac{D}{p}) + (\frac{D}{q}) = 0  \\
2 \quad \text{if} \ (\frac{D}{p}) = (\frac{D}{q}) = -1.   \\
\end{array} \right.  $  \\
Case for $  \delta _{a}. $ \ Since $ S_{a} = S_{ar} $ is given
above, we have \\
$  \delta _{a} = \sum _{v \in S_{a}} \delta _{v} =
\left \{\begin{array}{l} \delta _{2} \quad \text{if} \
D \equiv 3, 5, 7 (\text{mod} 8)  \\
0 \quad \text{if} \
D \equiv 1 (\text{mod} 8).
\end{array} \right.  $
So the remainder is to compute the local norm index $ \delta _{2} $
when $ D \equiv 3, 5, 7 (\text{mod} 8). $ So we assume now
$ D \equiv 3, 5, 7 (\text{mod} 8). $ By the Theorem 7.6 in [KT, p.332]
(see also [Q1, p.5054]),
$$ \delta _{2} = \log _{2} \left( \frac{c_{2} c_{D, 2}}{c_{w}} \left(
\frac{\parallel \Delta _{2} \Delta _{D, 2} d (K_{w} / \Q_{2})^{-6}
\parallel _{\Q_{2}}}{\parallel \Delta _{w} \parallel
_{K_{w}}}\right)^{1 / 12} \right). $$
By Lemma 2.1(1) above, we have $ c_{2} = c_{D, 2} = 2, \Delta _{D, 2}
= 64p^{2}q^{2}D^{6}. $ Also, by the results in [Q1, p.5058], we have
$ d (K_{w} / \Q_{2}) = \left \{\begin{array}{l} D \quad \text{if} \
D \equiv 5 (\text{mod} 8)  \\
4D \quad \text{if} \
D \equiv 3, 7 (\text{mod} 8).
\end{array} \right. $
From these discussion together with the results of $ c_{w} $ and
$ \Delta _{w} $ in Lemma 3.2 above, one can work out
$ \delta _{2} $ as follows. \\
If $ D \equiv 5 (\text{mod} 8), $ then $ \delta _{2} = 1; $ \\
If $ D \equiv 7 (\text{mod} 8), $ then $ \delta _{2} =
\left \{\begin{array}{l} 2 \quad \text{if} \
p \equiv 1 (\text{mod} 4)  \\
1 \quad \text{if} \
p \equiv 3 (\text{mod} 4);
\end{array} \right. $ \\
If $ D \equiv 3 (\text{mod} 8), $ then $ \delta _{2} =
\left \{\begin{array}{l} 1 \quad \text{if} \
p \equiv 1 (\text{mod} 4)  \\
2 \quad \text{if} \
p \equiv 3 (\text{mod} 4).
\end{array} \right. $ \\
Now our conclusion follows.  \quad $ \Box $
\par  \vskip 0.2 cm

Recall that $ \amalg\hskip-7pt\amalg(E/K) $ is the Shafarevich-Tate
group of $ E/K. $ We have the following
explicit parity relation between $ \text{rank}(E(K)) $
and $ \text{dim}_{2} \amalg\hskip-7pt\amalg(E/K)[2]. $
\par  \vskip 0.2 cm

{\bf Theorem 3.4.} \ Let $ E / \Q $ be the elliptic curve in (1.1),
$ \mu = \pm 1, $ and $ K = \Q (\sqrt{\mu D}) $ be the quadratic
number field with square-free integer $ D = D_{1} \cdots D_{n} $ as
in (1.2) above. Denote $ \mu _{0} = (1 - \mu ) / 2. $ Then we have \\
(1) \ $ \text{rank}(E(K)) \equiv \mu _{0} + \text{dim}_{2}
\amalg\hskip-7pt\amalg(E/K)[2] \ (\text{mod} 2) $ if one of the
following six hypotheses holds: \\
(1a) \ $ D \equiv \mu (\text{mod} 8) $ and $ (\frac{\mu D}{p}) =
(\frac{\mu D}{q}); $ \\
(1b) \ $ D \equiv 3 \mu (\text{mod} 8), \  p \equiv 3 (\text{mod} 4)
$ and $ (\frac{\mu D}{p}) = (\frac{\mu D}{q}); $ \\
(1c) \ $ D \equiv 3 \mu (\text{mod} 8), \  p \equiv 1
(\text{mod} 4) $ and $ (\frac{\mu D}{p}) + (\frac{\mu D}{q}) = 0; $ \\
(1d) \ $ D \equiv 5 \mu (\text{mod} 8) $ and
$ (\frac{\mu D}{p}) + (\frac{\mu D}{q}) = 0; $ \\
(1e) \ $ D \equiv 7 \mu (\text{mod} 8), \  p \equiv 1 (\text{mod} 4)
$ and $ (\frac{\mu D}{p}) = (\frac{\mu D}{q}); $ \\
(1f) \ $ D \equiv 7 \mu (\text{mod} 8), \  p \equiv 3
(\text{mod} 4) $ and $ (\frac{\mu D}{p}) + (\frac{\mu D}{q}) = 0. $ \\
(2) \ $ \text{rank}(E(K)) \equiv \mu _{0} + 1 + \text{dim}_{2}
\amalg\hskip-7pt\amalg(E/K)[2] \ (\text{mod} 2) $ if one of the
following six hypotheses holds: \\
(2a) \ $ D \equiv \mu (\text{mod} 8) $ and $ (\frac{\mu D}{p}) +
(\frac{\mu D}{q}) = 0; $ \\
(2b) \ $ D \equiv 3 \mu (\text{mod} 8), \  p \equiv 1 (\text{mod} 4)
$ and $ (\frac{\mu D}{p}) = (\frac{\mu D}{q}); $ \\
(2c) \ $ D \equiv 3 \mu (\text{mod} 8), \  p \equiv 3
(\text{mod} 4) $ and $ (\frac{\mu D}{p}) + (\frac{\mu D}{q}) = 0; $ \\
(2d) \ $ D \equiv 5 \mu (\text{mod} 8) $ and
$ (\frac{\mu D}{p}) = (\frac{\mu D}{q}); $ \\
(2e) \ $ D \equiv 7 \mu (\text{mod} 8), \  p \equiv 3 (\text{mod} 4)
$ and $ (\frac{\mu D}{p}) = (\frac{\mu D}{q}); $ \\
(2f) \ $ D \equiv 7 \mu (\text{mod} 8), \  p \equiv 1 (\text{mod} 4)
$ and $ (\frac{\mu D}{p}) + (\frac{\mu D}{q}) = 0. $
\par  \vskip 0.1 cm

{\bf Proof.} \ By Theorem 1 of [Kr, p.130], we have \\
$ \text{rank}(E(K)) \equiv \sum _{v \in M_{\Q}} \delta _{v}
+ \text{dim}_{2} \amalg\hskip-7pt\amalg(E/K)[2] =  \delta (E, \Q, K)
+ \text{dim}_{2} \amalg\hskip-7pt\amalg(E/K)[2] \ (\text{mod} 2). $
So the conclusion follows from Theorem 3.3 above.  \quad $ \Box $
\par  \vskip 0.2 cm

{\bf Corollary 3.5.} \ Let $ E / \Q $ and $ K $ be as in Theorem 3.4
above. If $ \sharp \amalg\hskip-7pt\amalg(E/K)[2] $ is a square
integer, then under one of the conditions in (2) for $ \mu = 1 $ (or
in (1) for $ \mu = -1 $) of Theorem 3.4, we have $ \text{rank}(E(K))
> 0. $
\par  \vskip 0.1 cm

{\bf Proof.} \ Obvious. \quad $ \Box $
\par  \vskip 0.2 cm

Now for an elliptic curve $ E $ over a number field $ F, $ and
a positive integer $ m, $ let $ \text{Sel}_{m}(E/F) $ be the
$ m-$Selmer group of $ E/F $ (see [Sil1, Chapt.10]).
\par  \vskip 0.2 cm

{\bf Corollary 3.6.} \ For the elliptic curves $ E / \Q $ in
(1.1) and $ E_{D} / \Q $ in (1.2) above, let $ \mu $ and $ \mu _{0} $
be as in Theorem 3.4 above. Then we have \\
(1) \ $ \text{dim}_{2} \text{Sel}_{2}(E_{\mu D}/\Q) \equiv \mu _{0}
+ \text{dim}_{2} \text{Sel}_{2}(E/\Q) \ (\text{mod} 2) $
if one of the six hypotheses in (1) of Theorem 3.4 above holds. \\
(2) \ $ \text{dim}_{2} \text{Sel}_{2}(E_{D}/\Q) \equiv \mu _{0} + 1
+ \text{dim}_{2} \text{Sel}_{2}(E/\Q) \ (\text{mod} 2) $ if one of
the six hypotheses in (2) of Theorem 3.4 above holds.
\par  \vskip 0.1 cm

{\bf Proof.} \ Let $ K =\Q (\sqrt{\mu D}) $ be as in Theorem 3.4
above. By Kramer's theorem (see [MR, Thm.2.7]), we have \\
$ \text{dim}_{2} \text{Sel}_{2}(E_{\mu D}/\Q) \equiv \text{dim}_{2}
\text{Sel}_{2}(E/\Q) + \delta (E, \Q, K) \ (\text{mod} 2). $ So the
conclusion follows from Theorem 3.3 above. \quad $ \Box $
\par  \vskip 0.2 cm

For an elliptic curve $ E / \Q, $ let $ L(E/ \Q, s) $ be its $ L-$
function (see [Sil1]). We denote its analytic rank by $ r_{an}(E/\Q), $
i.e., $ r_{an}(E/\Q) = \text{ord}_{s = 1} L(E/ \Q, s), $ which is the
order of $ L(E/ \Q, s) $ vanishing at $ s = 1. $
\par  \vskip 0.2 cm

{\bf Theorem 3.7.} \ Let $ E_{D} / \Q $ be the elliptic curve in (1.2)
above ($ E_{1} = E $ in (1.1) when take $ D = 1 $). Assume that one
of the following four hypotheses holds: \\
(1) \ $ p > 37 $ and the $ p-$Selmer group $ \text{Sel}_{p}(E_{D}/ \Q) $
is trivial; \\
(2) \ $ p > 37 $ and the $ q-$Selmer group $ \text{Sel}_{q}(E_{D}/ \Q) $
is trivial; \\
(3) \ $ 5 \nmid pqD, \ E_{D}[5] $ is an irreducible $ G_{\Q}-$module, and the
$ 5-$Selmer group $ \text{Sel}_{5}(E_{D}/ \Q) $ is trivial; \\
(4) \ $ 7 \nmid pqD, \ p \equiv 1, 4 \ (\text{mod} 7), \ E_{D}[7] $ is an
irreducible $ G_{\Q}-$module, and the $ 7-$Selmer group
$ \text{Sel}_{7}(E_{D}/ \Q) $ is trivial. \\
Then the rank and analytic rank of $ E_{D} / \Q $ are both equal to $ 0, $ i.e.,
$ \text{rank}(E_{D}(\Q)) =  r_{an}(E_{D}/\Q) = 0. $
\par  \vskip 0.1 cm

{\bf Proof.} \ First, assume (1) (resp. (2)), then \\
(a) \ $ E_{D} $ has multiplicative reduction at both $ p $ and $ q; $ \\
(b) \ Since $ E_{D} $ has no complex multiplication, by the work of [Ma1]
(or see[Cha, p.175]), for $ p > 37, $ both $ E_{D}[p] $ and $ E_{D}[q] $ are
irreducible $ G_{\Q}-$modules; \\
(c) \ By Prop.2.6 above, $ E_{D}[p] $ is ramified at $ q, $ and
$ E_{D}[q] $ is ramified at $ p; $ \\
(d) \ By assumption, $ \text{Sel}_{p}(E_{D}/ \Q) $ (resp.
$ \text{Sel}_{q}(E_{D}/ \Q)$ ) is trivial. \\
So all the conditions (a), (b), (c), (d) in Theorem 5 of [BSZ, p.3]
hold, and the conclusion follows. \\
Next,  assume (3) (resp. (4)), then \\
(a) \ By Cor.2.2 above, $ E_{D} $ has good ordinary reduction at $ 5 $
(resp. $ 7 $); \\
(b) \ $ E_{D}[5] $ (resp. $ E_{D}[7] $) is an irreducible $ G_{\Q}-$module; \\
(c) \ By Prop.2.6 above, $ E_{D}[5] $ (resp. $ E_{D}[7] $) is ramified at
$ p; $ \\
(d) \ $ \text{Sel}_{5}(E_{D}/ \Q) $ (resp. $ \text{Sel}_{7}(E_{D}/ \Q) $)
is trivial. \\
So all the conditions (a), (b), (c), (d) in Theorem 5 of [BSZ, p.3]
hold, and the conclusion follows.  \quad $ \Box $
\par  \vskip 0.2 cm

{\bf Theorem 3.8.} \ Let $ E_{D} / \Q $ be the elliptic curve in (1.2)
above ($ E_{1} = E $ in (1.1) when take $ D = 1 $). Assume that one
of the following two hypotheses holds: \\
(1) \ $ 5 \nmid pqD, \ E_{D}[5] $ is an irreducible $ G_{\Q}-$module, and the
$ 5-$Selmer group $ \text{Sel}_{5}(E_{D}/ \Q) $ has order $ 5; $ \\
(2) \ $ 7 \nmid pqD, \ p \equiv 1, 4 \ (\text{mod} 7), \ E_{D}[7] $ is an
irreducible $ G_{\Q}-$module, and the $ 7-$Selmer group
$ \text{Sel}_{7}(E_{D}/ \Q) $  has order $ 7. $ \\
Then the rank and analytic rank of $ E_{D} / \Q $ are both equal to $ 1, $ i.e.,
$ \text{rank}(E_{D}(\Q)) =  r_{an}(E_{D}/\Q) = 1. $
\par  \vskip 0.1 cm

{\bf Proof.} \ Assume (1) (resp. (2)), then \\
(a) \ By Cor.2.2 above, $ E_{D} $ has good ordinary reduction at $ 5 $
(resp. $ 7 $); \\
(b) \ $ E_{D}[5] $ (resp. $ E_{D}[7] $) is an irreducible $ G_{\Q}-$module; \\
(c) \ By Prop.2.6 above, $  E_{D}[5] $ (resp. $E_{D}[7] $) is ramified at
$ l $ for $ l = p $ or $ q; $ \\
(d) \ The conductor $ N $ of $ E_{D} $ is obviously not square-free,
and there are two distinct prime factors $ l \parallel N $
(i.e., $p, q $) such that $ E_{D}[5] $) (resp. $E_{D}[7]$)is ramified at
$ l; $  \\
(e) \ $ E_{D} $ obviously has good reduction at $ 5 $ (resp. $ 7 $); \\
(f) \ $ \text{Sel}_{5}(E_{D}/ \Q) $ (resp. $ \text{Sel}_{7}(E_{D}/ \Q) $)
has order $ 5 $ (resp. $ 7. $) \\
So all the conditions (a), (b), (c), (d), (e), (f) in Theorem 9 of
[BSZ, p.4] hold, and the conclusion follows.  \quad $ \Box $
\par  \vskip 0.2 cm

{\bf Remark.} \ For the elliptic curve $ E_{D} $ in Theorem 3.8 above,
since its conductor $ N = 2^{5}pqD^{2} $ has two distinct prime
factors of order one, i.e., $ p $ and $ q, $ by Theorem 1.5 of
[Zh, p.8], we know that the following two statements are equivalent: \\
(1) \ $  \text{rank}(E_{D}(\Q)) = 1 $ and
$ \sharp \amalg\hskip-7pt\amalg(E_{D}/ \Q) < + \infty ; $ \\
(2) \ $ r_{an}(E_{D}/\Q) = 1. $

\par     \vskip  0.4 cm

\hspace{-0.8cm}{\bf 4. Iwasawa theory for $ E_{D} $}

\par  \vskip 0.2 cm

Let $ E $ be an elliptic curve defined over a number field $ F, \
m $ be a positive integer and $ l $ be a prime number. Then for
any place $ v \in M_{F}, $ we have the Kummer homomorphisms \\
$ \kappa _{v, m} : \ E(F_{v}) \otimes \Z / m \Z \rightarrow
\text{H}^{1}(F_{v}, E[m]), $ \ and \ $
\kappa _{v, l^{\infty }} : \ E(F_{v}) \otimes \Q_{l} / \Z_{l}
\rightarrow \text{H}^{1}(F_{v}, E[l^{\infty }]), $ \\
where $ \Z_{l} $ is the ring of $ l-$adic integers and
$ E[l^{\infty }] $ is the $l-$primary torsion subgroup of $ E. $
Recall that the $ m-$Selmer group
$ \text{Sel}_{m}(E / F) $ of $ E / F $ is
defined as \\
$ \text{Sel}_{m}(E / F) = \text{ker} \{\text{H}^{1}(F, E[m])
\longrightarrow \prod _{v \in M_{F}} \text{H}^{1}(F_{v}, E[m])
/ \text{Im}(\kappa _{v, m}) \}, $ \\
and the $ l^{\infty }-$Selmer group
$ \text{Sel}_{l^{\infty }}(E / F) $ is defined as \\
$ \text{Sel}_{l^{\infty }}(E / F) =
\text{ker} \{\text{H}^{1}(F, E[l^{\infty }]) \longrightarrow
\prod _{v \in M_{F}} \text{H}^{1}(F_{v}, E[l^{\infty }])
/ \text{Im}(\kappa _{v, l^{\infty }}) \}. $ \\
Note that the $ l^{\infty }-$Selmer group can be defined for $ E $
over any algebraic extension $ M $ of $ \Q $ (see [Gr, p.63]).
There is a natural surjective homomorphism (see [Zh, p.3]) \
$$ \text{Sel}_{l}(E / F) \longrightarrow
\text{Sel}_{l^{\infty }}(E / F)[l], $$ and the properties of
$ \text{Sel}_{l^{\infty }}(E / F) $ can sometimes be deduced from
the ones of $ \text{Sel}_{l}(E / F) $ (see [BS, p.6]). \\
Let $ \Q_{\infty } $ be a $ \Z_{l}-$extension, i.e., it is a Galois
extension of $ \Q $ such that $ \text{Gal}(\Q_{\infty } / \Q) \cong
\Z_{l}, $ the additive group of $ l-$adic integers. So we have
$ \Q_{\infty } = \cup _{n \geq 0} \Q_{n}, $ where for each
$ n, \Q_{n} $ is a cyclic extension of $ \Q $ of degree $ l^{n} $
and $ \Q = \Q_{0} \subset \Q_{1} \subset \cdots \subset \Q_{n} \subset
\cdots . $ We write $ \Gamma = \text{Gal}(\Q_{\infty } / \Q), $ and
let $ \gamma \in \Gamma $ be a fixed topological generator. The completed
group ring $ \Lambda = \Z_{l}[[\Gamma ]] \cong \Z_{l}[[T]], $ where the
indeterminate $ T $ is identified with $ \gamma - 1. $ We write
$ \Gamma _{n} = \text{Gal}(\Q_{\infty } / \Q_{n}), $ then $ \Gamma _{n}
= \Gamma ^{l^{n}}. $ For the structure of the Iwasawa algebra $ \Lambda , $
see [Wa]. For an elliptic curve $ E $ defined over $ \Q, $ the Pontryagin
dual of its $ l^{\infty }-$Selmer group
$ \text{Sel}_{l^{\infty }}(E / \Q_{\infty }) $ is denoted by
$ X(E / \Q_{\infty }) =
\text{Hom}( \text{Sel}_{l^{\infty }}(E / \Q_{\infty }), \Q_{l}/\Z_{l}). $
It is a $ \Lambda-$module via the natural action of $ \Gamma $ on the group
$ \text{H}^{1}(\Q_{\infty }, E[l^{\infty }]), $ and one says that
$ \text{Sel}_{l^{\infty }}(E / \Q_{\infty }) $ is $ \Lambda-$cotorsion if
$ X(E / \Q_{\infty }) $ is $ \Lambda-$torsion (see [Gr, p.55]).
\par  \vskip 0.2 cm

Now let $ E_{D} / \Q $ be the elliptic curve in (1.2)
above ($ E_{1} = E $ in (1.1) when take $ D = 1 $). Assume that the
prime number $ l $ satisfies one of the following two hypotheses: \\
(1) \ $ l = 5 $ and $ 5 \nmid pqD; $ \\
(2) \ $ l = 7, \ 7 \nmid pqD, $ and $ p \equiv 1, 4 (\text{mod} 7). $ \\
Then by Cor.2.2 above, $ E_{D} $ has good ordinary reduction at such $ l. $
So by Mazur's control theorem (see [Gr, p.54]), the natural maps
$$ \text{Sel}_{l^{\infty }}(E_{D} / \Q_{n})
\longrightarrow \text{Sel}_{l^{\infty }}(E_{D} / \Q_{\infty })^{\Gamma _{n}} $$
have finite kernel and cokernel, of bounded order as $ n $ varies. \\
Such $ E_{D} / \Q $ also has multiplicative reduction at $ p $ and $ q, $
so for the prime number $ l $ such that $ l = p, q $ or satisfies one of the
above two hypotheses (1) and (2), by Kato-Rohrlich's theorem (see [Gr, p.55]),
we know that $ \text{Sel}_{l^{\infty }}(E_{D} / \Q_{\infty }) $ is
$ \Lambda-$cotorsion. \\
Furthermore, under this hypothesis, we have the following results.
\par  \vskip 0.2 cm

{\bf Proposition 4.1.} \ Let $ E_{D} / \Q $ be the elliptic curve in (1.2)
above ($ E_{1} = E $ in (1.1) when take $ D = 1 $). Let $ l $ be a
prime number satisfying one of the following two hypotheses: \\
(1) \ $ l = 5 $ and $ 5 \nmid pqD; $ \\
(2) \ $ l = 7, \ 7 \nmid pqD, $ and $ p \equiv 1, 4 (\text{mod} 7). $ \\
Then the map
$$ \text{Sel}_{l^{\infty }}(E_{D} / \Q)
\longrightarrow \text{Sel}_{l^{\infty }}(E_{D} / \Q_{\infty })^{\Gamma } $$
is surjective. If $ \text{Sel}_{l^{\infty }}(E_{D} / \Q) = 0, $
then $ \text{Sel}_{l^{\infty }}(E_{D} / \Q_{\infty }) = 0 $
also.
\par  \vskip 0.1 cm

{\bf Proof.} \ By Cor.2.2 above, $ E_{D} $ has good ordinary reduction at
such $ l; $ by Lemma 2.3 above, we have
$ l \nmid \sharp \widetilde{E_{D,}}_{l}(\F_{l}); $ and by Lemma 2.1,
$ l \nmid c_{l^{\prime }} $ for any prime number $ l^{\prime }. $ So
the conditions (i), (ii), (iii) of Prop.3.8 in [Gr, p.80] hold (see
also the Remark there), and the conclusion follows. \quad $ \Box $
\par  \vskip 0.2 cm

{\bf Proposition 4.2.} \ Let $ E_{D} / \Q $ be the elliptic curve in (1.2)
above ($ E_{1} = E $ in (1.1) when take $ D = 1 $). Let $ l $ be a
prime number satisfying one of the following three hypotheses: \\
(1) \ $ l = p $ or $ q; $ \\
(2) \ $ l = 5 $ and $ 5 \nmid pqD; $ \\
(3) \ $ l = 7, \ 7 \nmid pqD, $ and $ p \equiv 1, 4 (\text{mod} 7). $ \\
Then for all $ n \geq 0, $ the map \
$ \text{Sel}_{l^{\infty }}(E_{D} / \Q_{n})
\longrightarrow \text{Sel}_{l^{\infty }}(E_{D} / \Q_{\infty }) $
is injective. Moreover,
$$ \text{corank}_{\Z_{l}}(\text{Sel}_{l^{\infty }}(E_{D} / \Q_{\infty }))
\equiv \text{corank}_{\Z_{l}}(\text{Sel}_{l^{\infty }}(E_{D} / \Q))
(\text{mod} 2). $$
\par  \vskip 0.1 cm

{\bf Proof.} \ Under our assumption, $ E_{D} $ has good ordinary
or multiplicative reduction at $ l. $ Also, by the above discussion,
we know that $ \text{Sel}_{l^{\infty }}(E_{D} / \Q_{\infty }) $ is
$ \Lambda-$cotorsion, so the conclusion follows from the Prop.3.9
and Prop.3.10 of [Gr, pp.81, 82]. \quad $ \Box $
\par  \vskip 0.2 cm

Now for the elliptic curves $ E_{D} / \Q $ and the prime number $ l $
as in the above Proposition 4.2, by Mazur and Swinnerton-Dyer's construction,
there is an element
$ \mathfrak{L}(E_{D} / \Q, T) \in \Lambda \otimes _{\Z_{l}} \Q_{l} $
with some interpolation property, from which one can define the
$ l-$adic $L-$ function $ L_{l}(E_{D} / \Q, s). $ For the general
theory of $ l-$adic $ L-$function of elliptic curves, see [MSD] and [Gr].
By Weierstrass' preparation theorem, we have
$ \mathfrak{L}(E_{D} / \Q, T) = l^{m_{1}} \cdot U(T) \cdot f(T), $ where
$ f(T) $ is a distinguished polynomial, $ U(T) $ is an invertible power
series and $ m_{1} \in \Z. $ As in [GV, pp.19, 20], we write
$ f_{E_{D}}^{\text{anal}}(T) = l^{m_{1}} \cdot f(T). $ On the other hand, since
$ \text{Sel}_{l^{\infty }}(E_{D} / \Q_{\infty }) $ is
$ \Lambda-$cotorsion, i.e., $ X(E_{D} / \Q_{\infty }) $ is
$ \Lambda-$torsion, one has a pseudo-isomorphism
$$ X(E_{D} / \Q_{\infty }) \sim
(\oplus _{i = 1}^{n} \Lambda / (f_{i}(T)^{a_{i}})) \oplus
(\oplus _{j = 1}^{m} \Lambda / (l^{b_{j}})), $$
where $ f_{i}(T) $ are irreducible distinguished polynomials in
$ \Lambda , $ and $ a_{i}, b_{j} $ are non-negative integers.
Then the characteristic polynomial for the
$ \Lambda-$module $ X(E_{D} / \Q_{\infty }) $
is defined by $ f_{E_{D}}^{\text{alg}}(T) =
l^{m_{2}} \cdot \prod _{i = 1}^{n} f_{i}(T)^{a_{i}}, $ where
$ m_{2} = \sum _{j = 1}^{m}b_{j}. $ By Kato's theorem about the
main conjecture (see [GV, p.21]), the polynomial
$ f_{E_{D}}^{\text{alg}}(T) $ divides $ f_{E_{D}}^{\text{anal}}(T) $
in $ \Q_{l}[T]. $ Moreover, by Greenberg's theorem (see [Gr, p.61]),
the characteristic ideal of $ X(E_{D} / \Q_{\infty }) $
is fixed by the involution $ \iota $ of $ \Lambda $ induced by
$ \iota (\sigma ) = \sigma ^{-1} $ for all $ \sigma \in \Gamma . $
\par  \vskip 0.2 cm

{\bf Theorem 4.3.} \ Let $ E_{D} / \Q $ be the elliptic curve in (1.2)
above ($ E_{1} = E $ in (1.1) when take $ D = 1 $). Let $ l $ be a
prime number satisfying one of the following three hypotheses: \\
(1) \ $ l = p $ or $ q; $ \\
(2) \ $ l = 5 $ and $ 5 \nmid pqD; $ \\
(3) \ $ l = 7, \ 7 \nmid pqD, $ and $ p \equiv 1, 4 (\text{mod} 7). $ \\
Then $ \text{Sel}_{l^{\infty }}(E_{D} / \Q_{\infty }) $ has no proper
$ \Lambda-$submodules of finite index. In particular, if
$ \text{Sel}_{l^{\infty }}(E_{D} / \Q_{\infty }) \neq 0, $ then
$ \text{Sel}_{l^{\infty }}(E_{D} / \Q_{\infty }) $ is finite. \\
Moreover, for $ l $ satisfying the hypothesis (2) or (3) here, if
$ \text{Sel}_{l^{\infty }}(E_{D} / \Q) $
is finite, then $ f_{E_{D}}^{\text{alg}}(0) \sim \sharp
\text{Sel}_{l^{\infty }}(E_{D} / \Q). $ Here, for
$ a, b \in \Q_{l}^{\times }, $ we write $ a \sim b $ to indicate that
$ a $ and $ b $ have the same $ l-$adic valuation.
\par  \vskip 0.1 cm

{\bf Proof.} \ By Lemma 2.1(5) above, the torsion subgroup
$ E_{D}(\Q)_{\text{tors}} \cong \Z / 2 \Z \times \Z / 2 \Z, $
so for the prime number $ l $ under our assumption,
$ E_{D}(\Q)_{\text{tors}}[l^{\infty }] = 0. $ Also, by the above
discussion, we know that
$ \text{Sel}_{l^{\infty }}(E_{D} / \Q_{\infty }) $ is
$ \Lambda-$cotorsion, so our first conclusion follows from the
Prop.4.14 of [Gr, p.102]. \\
Next we come to show our second conclusion. As
$ \text{Sel}_{l^{\infty }}(E_{D} / \Q_{\infty }) $ is
$ \Lambda-$cotorsion, let $  f_{E_{D}}^{\text{alg}}(T) $ be its
characteristic polynomial as above, i.e., $ f_{E_{D}}^{\text{alg}}(T) $
is a generator of the characteristic ideal of the $ \Lambda-$module
$ X(E_{D} / \Q_{\infty }), $ the Pontryagin dual of
$ \text{Sel}_{l^{\infty }}(E_{D} / \Q_{\infty }). $ Denote
$ \theta _{n} = \gamma ^{l^{n}} - 1 = (1 + T)^{l^{n}} - 1 \in \Lambda $
for each $ n \geq 0. $ We know,
$ X(E_{D} / \Q_{\infty }) / \theta _{n} X(E_{D} / \Q_{\infty }) $ is the
Pontryagin dual of
$ \text{Sel}_{l^{\infty }}(E_{D} / \Q_{\infty })^{\Gamma _{n}}, $ and the
torsion subgroup of
$ X(E_{D} / \Q_{\infty }) / \theta _{n} X(E_{D} / \Q_{\infty }) $ is then
dual to \\
$ \text{Sel}_{l^{\infty }}(E_{D} / \Q_{\infty })^{\Gamma _{n}} /
(\text{Sel}_{l^{\infty }}(E_{D} / \Q_{\infty })^{\Gamma _{n}})_{\text{div}} $
(see [Gr, p.82]), \\
In particular,
$ X(E_{D} / \Q_{\infty }) / T X(E_{D} / \Q_{\infty }) $ is the
Pontryagin dual of
$ \text{Sel}_{l^{\infty }}(E_{D} / \Q_{\infty })^{\Gamma }. $ As assumed,
$ \text{Sel}_{l^{\infty }}(E_{D} / \Q) $ is finite, and so by the
above discussion,
$ \text{Sel}_{l^{\infty }}(E_{D} / \Q_{\infty })^{\Gamma } $ is also
finite, hence $ X(E_{D} / \Q_{\infty }) / T X(E_{D} / \Q_{\infty }) $
is finite. Therefore, $ T \nmid f_{E_{D}}^{\text{alg}}(T), $ so
$ f_{E_{D}}^{\text{alg}}(0) \neq 0. $ In the following, For an element
$ c \in \Z_{l}, $ the highest power of $ l $ dividing $ c $ is denoted
by $ c^{(l)}. $ \\
Now we assume that $ l $ satisfies the hypothesis (2), i.e., $ l = 5 $ and
$ 5 \nmid pqD. $ Then $ E_{D} $ has good ordinary reduction at $ 5, $ and by
Lemma 2.3 above, $  \sharp \widetilde{E_{D,}}_{5}(\F_{5}) = 4 $ or $ 8. $
So $ \widetilde{E_{D,}}_{5}(\F_{5})[5^{\infty }] = 0. $ Also by
Lemma 2.1, we have $ c_{l^{\prime }} = 2 $ or $ 4 $ for any
$ l^{\prime } \mid N_{E_{D}}, $ the conductor of $ E_{D}, $ and
$ E_{D}(\Q)_{\text{tors}} \cong \Z / 2 \Z \times \Z / 2 \Z. $ So
$ c_{l^{\prime }}^{(5)} = 1 $ for any $ l^{\prime } \mid N, $ and
$ E_{D}(\Q)[5^{\infty }] = 0. $ Hence by Theorem 4.1 of [Gr, p.85], we
get  \begin{align*} f_{E_{D}}^{\text{alg}}(0)
&\sim (\prod _{l^{\prime } \mid N_{E_{D}}} c_{l^{\prime }}^{(5)}) \cdot
(\sharp \widetilde{E_{D,}}_{5}(\F_{5})[5^{\infty }])^{2} \cdot
\sharp \text{Sel}_{5^{\infty }}(E_{D} / \Q) / ( \sharp
E_{D}(\Q)[5^{\infty }])^{2} \\
&= 1 \cdot 1^{2} \cdot \sharp \text{Sel}_{5^{\infty }}(E_{D} / \Q) /
1^{2} = \sharp \text{Sel}_{5^{\infty }}(E_{D} / \Q),
\end{align*}
i.e., $ f_{E_{D}}^{\text{alg}}(0) \sim
\sharp \text{Sel}_{5^{\infty }}(E_{D} / \Q). $
The case for $ l $ satisfying the hypothesis (3) can
be similarly done, and the proof is completed.
\quad $ \Box $
\par  \vskip 0.1 cm

{\bf Remark.} \ For the elliptic curve $ E_{D} / \Q $ in (1.2)
above, for every prime number $ l > 2, $ by Lemma 2.1 above,
we have $ E_{D}(\Q)[l^{\infty }] = 0, $ so
$ E_{D}(\Q_{\infty })[l^{\infty }] = 0 $ because $ \Gamma $ is
pro-$l $ (see [Gr, p.102, line -10]). so
$ E_{D}(\Q_{\infty })_{\text{tors}} $ is a $ 2-$group, i.e.,
its every element is of $ 2-$power order.
\par  \vskip 0.2 cm

For the elliptic curve $ E_{D} / \Q $ as in (1.2) above, let
$ \Omega _{D} $ be its N\'{e}ron period. Now we let $ l $ be a
prime number satisfying one of the following two hypotheses: \\
(1) \ $ l = 3 $ and $ 3 \nmid pqD; $ \\
(2) \ $ l = 7, \ 7 \nmid pqD, $ and $ p \equiv 2, 3, 6 (\text{mod} 7). $
\par  \vskip 0.2 cm

Then by Cor.2.2 above, we know that $ E_{D}$ has good supersingular
reduction at such $ l. $ By Lemma 2.1 above, we have $ c_{l^{\prime }} =
2 $ or $ 4 $ for any prime number $ l^{\prime} \mid N_{E_{D}} = 2^{5}pqD^{2}, $
so our $ l \nmid \text{Tam}(E_{D} / \Q ) = \prod_{l^{\prime} < \infty }
c_{l^{\prime }}. $ Also by Prop.2.7 above, we have $ \rho _{l}(G_{\Q}) =
\text{Gl}_{2}(\F_{l}). $ Therefore, if
$ \text{ord}_{l} (L(E_{D} / \Q, 1) / \Omega _{D}) = 0, $ then over
the $ \Z_{l}-$extension $ \Q_{\infty } / \Q $ as above, by The
Theorem 0.1 of [Ku, p.196], we have the following conclusion: \\
(1) \ $ (\amalg\hskip-7pt\amalg(E_{D}/ \Q_{\infty })[l^{\infty }])^{\wedge }
\cong \Lambda  $ as $ \Lambda-$modules, where
$ (\amalg\hskip-7pt\amalg(E_{D}/ \Q_{\infty })[l^{\infty }])^{\wedge } $
is the Pontryagin dual of
$ \amalg\hskip-7pt\amalg(E_{D}/ \Q_{\infty })[l^{\infty }]; $ \\
(2) \ $ \text{rank}(E_{D}( \Q_{n})) = 0 $ and
$ \sharp \amalg\hskip-7pt\amalg(E_{D}/ \Q_{n})[l^{\infty }] =
l^{e_{n}} $ with $ e_{n} = [\frac{l^{n + 1}}{l^{2} - 1} - \frac{n}{2}] $
for any $ n \geq 0; $ \\
(3) \ $ (\amalg\hskip-7pt\amalg(E_{D}/ \Q_{n})[l^{\infty }])^{\wedge }
\cong \Z_{l}[\text{Gal}(\Q_{n} / \Q)] / (\theta _{\Q_{n}},
v_{n-1, n}(\theta _{\Q_{n-1}})) $ as \\
$ \Z_{l}[\text{Gal}(\Q_{n} / \Q)]-$modules
for any $ n \geq 0, $ where $ \theta _{\Q_{n}} $ is the modular element of
Mazur and Tate (see [Ku] for the detail). \\
In fact, the Mordell-Weil group $ E_{D}( \Q_{n}) $ in the above result (2)
can be determined as follows.
\par  \vskip 0.2 cm

{\bf Theorem 4.4.} \ Let $ E_{D} / \Q $ be the elliptic curve in (1.2)
above ($ E_{1} = E $ in (1.1) when take $ D = 1 $). Let $ l $ be a
prime number satisfying one of the following two hypotheses: \\
(1) \ $ l = 3 $ and $ 3 \nmid pqD; $ \\
(2) \ $ l = 7, \ 7 \nmid pqD, $ and $ p \equiv 2, 3, 6 (\text{mod} 7). $ \\
If $ \text{ord}_{l} (L(E_{D} / \Q, 1) / \Omega _{D}) = 0, $ then over
the $ \Z_{l}-$extension $ \Q_{\infty } / \Q $ as above, we have
$  E_{D} (\Q_{n}) \cong \Z / 2 \Z \times Z / 2 \Z $ for all $ n \geq 0. $
\par  \vskip 0.1 cm

{\bf Proof.} \ By the above discussion, we know that
$ \text{rank}(E_{D}( \Q_{n})) = 0. $ So $ E_{D} (\Q_{n}) =
E_{D} (\Q_{n})_{\text{tors}}. $ Since $ E_{D} $ has good supersingular
reduction at such $ l, \ E_{D}(\Q(\mu _{l^{n + 1}})) $ does not contain
a point of order $ l $ for any $ n \geq 0 $ (see [Ku, p.200, line-2]),
where $ \mu _{l^{n + 1}} $ is the group of $ l^{n + 1}-$th roots of
unity. Since $ \Q_{\infty } $ is in fact the cyclotomic $ \Z_{l}-$extension
of $ \Q, $ we have $ \Q_{n} \subset \Q(\mu _{l^{n + 1}}), $ and so
$ E_{D} (\Q_{n})[l^{\infty }] = 0 $ for any $ n \geq 0. $ On the other hand,
$ l $ is totally ramified in $ \Q_{n}. $ Let $ \mathfrak{p}_{n} $ be the
unique prime ideal of $ \Q_{n} $ lying over $ l, $ then the residue degree
$ f(\mathfrak{p}_{n} / l) = 1, $ and the residue field
$ k_{ \mathfrak{p}_{n}} = \F_{l}. $ So if $ l = 3, $ then by
Lemma 2.1(6) above, we have
$  E_{D} (\Q_{n})_{\text{tors}} / E_{D} (\Q_{n})[3^{\infty }]
\cong \Z / 2 \Z \times Z / 2 \Z, $ and then our conclusion follows because
$ E_{D} (\Q_{n})[3^{\infty }] = 0. $ If $ l = 7, $ then by Lemma 4.2(1)
of [QZ1, p.1379], we have $ \sharp E_{D} (\Q_{n})_{\text{tors}} \mid
\sharp \widetilde{E_{D,}}_{\mathfrak{p}_{n}}(\F_{7}) \cdot 7^{2 t_{7}} $
for some $ t_{7} \in \Z_{\geq 0}. $ By Lemma 2.3 above,
$ \sharp \widetilde{E_{D,}}_{\mathfrak{p}_{n}}(\F_{7}) = 8. $ Also,
by the above discussion, $ 7 \nmid \sharp E_{D} (\Q_{n})_{\text{tors}}. $
So $ \sharp E_{D} (\Q_{n})_{\text{tors}} \mid 8. $ Obviously,
$ E_{D} (\Q_{n})_{\text{tors}} \supset E_{D} (\Q_{n})[2] \cong
\Z / 2 \Z \times Z / 2 \Z, $ so $ E_{D} (\Q_{n})_{\text{tors}} \cong
\Z / 2 \Z \times Z / 2 \Z $ or $ \Z / 2 \Z \times Z / 4 \Z. $
The remainder is to show that $ E_{D} (\Q_{n})_{\text{tors}}
\ncong \Z / 2 \Z \times Z / 4 \Z, $ and this follows from the
following \\
Assertion. \ $ E_{D}(\Q(\mu _{7^{n}})) $ does not contain a point of
order $ 4 $ for any $ n \geq 0. $ \\
To see this, firstly, by Lemma 2.1 above, $ E_{D}(\Q)_{\text{tors}}
\cong \Z / 2 \Z \times Z / 2 \Z, $ so we may as well assume that
$ n > 0. $ Obviously $ E_{D}[2] = \{O, (0, 0), (-\varepsilon pD, 0),
(-\varepsilon qD, 0) \}, $ so $ E_{D}(\Q(\mu _{7^{n}})) $ contains a
point $ P_{4}$ of order $ 4 $ if and only if $ 2 P_{4} = (0, 0),
(-\varepsilon pD, 0) $ or $ (-\varepsilon qD, 0). $ And by Theorem 4.2
of [Kn, p.85], this is equivalent to say that
(we write $ F = \Q(\mu _{7^{n}}) $): \
(a) $ \varepsilon pD, \varepsilon qD \in F^{2}; $ or (b)
$ -\varepsilon pD, 2 \varepsilon D \in F^{2}; $ or (c)
$ -\varepsilon qD, - 2 \varepsilon D \in F^{2}. $ But all of these
cases are impossible because $ 7 $ is the unique prime number which
ramifies in $ F $ and $ 7 \nmid pq. $ So the above Assertion follows,
and the proof is completed.  \quad $ \Box $

\par     \vskip  0.4 cm

\hspace{-0.8cm}{\bf 5. $ L-$function, root number and
parity conjecture}

\par  \vskip 0.2 cm

Let $ E / \Q $ be the elliptic curve in (1.1), and its quadratic
$ D-$twist $ E_{D} / \Q $ in (1.2) above. Let $ K = \Q (\sqrt{D}) $
and $ K^{\prime } = \Q (\sqrt{-D}). $ The $ (-D)-$twist of such $ E $ is
$$  E_{-D}= E^{\varepsilon }_{-D}:  \ y^{2} = x(x - \varepsilon p D)(x -
\varepsilon q D). \quad \eqno(5.1) $$ So, $ E^{\varepsilon }_{-D} =
E^{- \varepsilon }_{D}. $ \\
As before, Let $ L(E/ \Q, s), \ L(E_{D}/ \Q, s) $ and $ L(E_{-D}/ \Q, s) $
be the $ L-$functions of $ E/ \Q, \ E_{D}/ \Q $ and $ E_{-D}/ \Q $
respectively, and write
\begin{align*} &L(E/ \Q, s)= \Sigma _{n = 1}^{\infty } a_{1}(n)n^{-s}, \
L(E_{D}/ \Q, s)= \Sigma _{n = 1}^{\infty } a_{D}(n)n^{-s}, \\
&L(E_{-D}/ \Q, s)= \Sigma _{n = 1}^{\infty } a_{-D}(n)n^{-s}
\end{align*}
with coefficients $ a_{1}(n), a_{D}(n), a_{-D}(n) $ respectively. Let
\begin{align*} &\Lambda (E/ \Q, s)=(\frac{\sqrt{N_{E}}}{2\pi })^{s}
\Gamma (s) L(E/ \Q, s), \ \Lambda (E_{D}/ \Q, s)=
(\frac{\sqrt{N_{E_{D}}}}{2\pi })^{s} \Gamma (s) L(E_{D}/ \Q, s), \\
&\Lambda (E_{-D}/ \Q, s)=
(\frac{\sqrt{N_{E_{-D}}}}{2\pi })^{s} \Gamma (s) L(E_{-D}/ \Q, s),
\end{align*}
where $ N_{E}, N_{E_{D}} $ and $ N_{E_{-D}} $ are the conductors of
$ E, E_{D} $ and $ E_{-D}, $ respectively. Since these curves are
modular over $ \Q, $ their $ L-$functions have analytic continuation
to $ \C $ and satisfy functional equations (see [Sil1, p.362]):
\begin{align*} &\Lambda (E/ \Q, 2 - s) = \omega _{E} \Lambda (E/ \Q, s), \
\Lambda (E_{D}/ \Q, 2 - s) = \omega _{E_{D}} \Lambda (E_{D}/ \Q, s), \\
&\Lambda (E_{-D}/ \Q, 2 - s) =
\omega _{E_{-D}} \Lambda (E_{-D}/ \Q, s),
\end{align*}
where $ \omega _{E}, \omega _{E_{D}}, \omega _{E_{-D}} \in \{1, -1 \} $
are the corresponding root numbers. Let
$ \chi _{K} $ and $ \chi _{K^{\prime }} $ be the quadratic Dirichlet
characters associated to $ K $ and $ K^{\prime }, $ respectively.
Then if $ ( d(K), 2 N_{E}) = 1, $ we have
$ L(E_{D}/ \Q, s) = L(E/ \Q, \chi _{K}, s) $
(see, e.g., [Kol1, p.524], [Kol2, p.475]). So $ L(E/K, s) = L(E / \Q, s)
\cdot L(E/ \Q, \chi _{K}, s) = L(E / \Q, s) \cdot
L(E_{D}/ \Q, s) $ (see also [DFK, p.186]), from which their root
numbers satisfy
$ \omega _{E/K} = \omega _{E/\Q} \cdot \omega _{E_{D}/\Q}. $
Similar for $ L(E_{-D}/ \Q, s). $ We write \\
$ L(E/ \Q, \chi _{K}, s) = \Sigma _{n = 1}^{\infty } a_{1}(n)
\chi _{K}(n) n^{-s} $ with coefficients $ a_{1}(n) \chi _{K}(n). $
\par  \vskip 0.2 cm

{\bf Lemma 5.1.} \ Assume that $ (D, 2pq) = 1. $ Then for the above
root numbers $ \omega _{E}, \omega _{E_{D}} $ and $ \omega _{E_{-D}}, $
we have \\
(1) \ if $ D \equiv 1 (\text{mod} 4), $ then $ \omega _{E_{D}} =
\chi _{K}(-2pq) \omega _{E}. $  \\
(2) \ if $ D \equiv 3 (\text{mod} 4), $ then $ \omega _{E_{-D}} =
\chi _{K^{\prime }}(-2pq) \omega _{E}. $
\par  \vskip 0.1 cm

{\bf Proof.} \ The discriminants of the quadratic number fields $ K $
and $ K^{\prime } $ are
$$ d(K) = \left \{\begin{array}{l} D \quad \text{if}
\ D \equiv 1 (\text{mod} 4) \\
4D \quad \text{if} \ D \equiv 3 (\text{mod} 4),
\end{array} \right. \ \text{and} \ d(K^{\prime }) =
\left \{\begin{array}{l} -4 D \quad \text{if}
\ D \equiv 1 (\text{mod} 4) \\
-D \quad \text{if} \ D \equiv 3 (\text{mod} 4),
\end{array} \right. $$
respectively. If $ ( d(K), N_{E}) = 1, $ then $ \omega _{E_{D}} =
\chi _{K}(- N_{E}) \omega _{E}, $ and if
$ ( d(K^{\prime }), N_{E}) = 1, $ then $ \omega _{E_{-D}} =
\chi _{K^{\prime }}(- N_{E}) \omega _{E} $ (see [DFK, p.186]).
Note that $  N_{E} = 2^{5}pq, $
the conclusion follows. \quad $ \Box $
\par \vskip 0.2 cm

The curve $ E / \Q $ in (1.1) above is $ 2-$isogeny to the following
elliptic curve
$$  E^{\prime } : \ y^{2} = x^{3} - 2 \varepsilon (p + q) x^{2} +
4 x, \quad \eqno(5.2) $$
and the isogeny is as follows. \\
$ \varphi : \ E \longrightarrow  E^{\prime }, \ (x, y) \mapsto
(x + \varepsilon (p + q) + pq \cdot x^{-1}, \ y - pqy \cdot x^{-2}). $ \\
This will be used in the following calculation of the root numbers.
Obviously, the conductor of $ E^{\prime } / \Q $ is $ N_{E\prime } =
N_{E} = 2^{5}pq, $ and the discriminant is $ \Delta _{E^{\prime }}
= 2^{12} pq. $ Firstly, we need the following result.
\par \vskip 0.2 cm

{\bf Lemma 5.2.} \ Let $ E^{\prime } / \Q $ be the elliptic curve in
(5.2) above. \\
(1) \ At each prime $ l \mid N_{E^{\prime}}, $ the Kodaira type is
as follows: \\
$ I_{3}^{\ast } $ for $ l = 2, $ and $ I_{1} $ for $l = p $ or $ q. $ \\
(2) \ The Tamagawa number $ c_{2} = 2 $ or $ 4, $ more precisely, \\
 $ c_{2} = 2 $ if one of the following three hypotheses holds: \\
(a) \ $ p \equiv 3 (\text{mod} 8); \ $ (b) \ $ \varepsilon = 1 $
and $ p \equiv 1 (\text{mod} 8); \ $ (c) \ $ \varepsilon = -1 $
and $ p \equiv 5 (\text{mod} 8). $ \\
$ c_{2} = 4 $ if one of the following three hypotheses holds: \\
(a$^{\prime }$) \ $ p \equiv 7 (\text{mod} 8); \ $ (b$^{\prime }$) \
$ \varepsilon = 1 $ and $ p \equiv 5 (\text{mod} 8); \ $ (c$^{\prime }$) \
$ \varepsilon = -1 $ and $ p \equiv 1 (\text{mod} 8). $ \\
(3) \ The Tamagawa numbers $ c_{p} = c_{q} = 1. $
\par \vskip 0.1 cm

{\bf Proof.} This is a consequence of direct calculation by
the Algorithm of [Ta]. \quad $ \Box $
\par \vskip 0.2 cm

Now we come to calculate the root numbers.
\par \vskip 0.2 cm

{\bf Theorem 5.3.} \ Let $ \omega _{E} $ be the root number of the
the elliptic curve $ E / \Q $ in (1.1) above. \\
(1) \ If $ \varepsilon = 1, $ then $ \omega _{E} =
\left \{\begin{array}{l} 1 \quad \text{if}
\ p \equiv 5, 7 \ (\text{mod} \ 8) \\
-1 \quad \text{if} \ p \equiv 1, 3 \ (\text{mod} \ 8);
\end{array} \right. $ \\
(2) \ If $ \varepsilon = -1, $ then $ \omega _{E} =
\left \{\begin{array}{l} 1 \quad \text{if}
\ p \equiv 3, 5 \ (\text{mod} \ 8) \\
-1 \quad \text{if} \ p \equiv 1, 7 \ (\text{mod} \ 8).
\end{array} \right. $
\par \vskip 0.1 cm

{\bf Proof.} \ To begin with, from [Roh, p.122], we have
$  \omega _{E} = \prod _{l \leq \infty } \omega _{l}, $ where
$ \omega _{l} = \pm 1 $ is the local root number. And by Prop.1
in [Roh1, p.123] one has $ \omega _{\infty } = -1, $ so
$  \omega _{E} = - \prod _{l < \infty } \omega _{l}. $ Since the
conductor is $ N_{E} = 2^{5} pq, $ for any prime number $ l \neq
2, p, q, E $ has good reduction at $ l, $ so by Prop.2(iv) in
[Roh, p.126], we have $ \omega _{l} = 1 $ for every such $ l. $
Also, since $ E / \Q $ has multiplicative reduction at both $ p $
and $ q, $ by discussion in Lemma 2.1 above, and by Prop.3(iii) in
[Roh, p.132], we have \\
(1) \ $ \omega _{p} = \omega _{q} = 1 $ if $ \varepsilon = 1 $ and
$ p \equiv 3, 5 \ (\text{mod} \ 8); $ \\
(2) \ $ \omega _{p} = \omega _{q} = -1 $ if $ \varepsilon = 1 $ and
$ p \equiv 1, 7 \ (\text{mod} \ 8); $ \\
(3) \ $ \omega _{p} = -1, \ \omega _{q} = 1 $ if $ \varepsilon = -1 $
and $ p \equiv 1, 3 \ (\text{mod} \ 8); $ \\
(4) \ $ \omega _{p} = 1, \ \omega _{q} = -1 $ if $ \varepsilon = -1 $
and $ p \equiv 5, 7 \ (\text{mod} \ 8). $ \\
So the remainder is the most difficult factor $ \omega _{2}. $ To
work out $ \omega _{2}, $ from [D], one can obtain the following formula
$$ \omega _{2} = \sigma _{\varphi }(E / \Q_{2}) \cdot
(\varepsilon (p + q), -pq)_{\Q_{2}}
\cdot (-2 \varepsilon (p + q), 4)_{\Q_{2}}, $$ recall that
$ ( , )_{\Q_{2}} $ is the Hilbert symbol (see [Se2, p.206]), $ \varphi $
is the isogeny in (5.2) above, and here,
$$ \sigma _{\varphi }(E / \Q_{2}) =
(-1)^{\text{ord}_{2}(\frac{\sharp \text{coker}\varphi _{2}}{\sharp
\text{ker}\varphi _{2}})} = (-1)^{1 + \text{ord}_{2} \sharp
\text{coker}\varphi _{2}}, $$ where $ \varphi _{2} : \ E(\Q_{2})
\longrightarrow E^{\prime }(\Q_{2}) $ is the local homomorphism induced
by $ \varphi . $ Since $ ( , )_{\Q_{2}} $ is biadditive, we have
$ (-2 \varepsilon (p + q), 4)_{\Q_{2}} =
(-2 \varepsilon (p + q), 2)_{\Q_{2}}^{2} = 1, $ so
$ \omega _{2} = \sigma _{\varphi }(E / \Q_{2}) \cdot
(\varepsilon (p + q), -pq)_{\Q_{2}}. $ To calculate
$ (\varepsilon (p + q), -pq)_{\Q_{2}}, $ we consider the equation
$ \varepsilon (p + q) x^{2} - pq y^{2} = 1. $ Let $ f(x, y) =
\varepsilon (p + q) x^{2} - pq y^{2} - 1, $ then
$ \frac{\partial f}{\partial y}(x, y) = -2pqy, $ and it is easy to
see that $ \text{ord}_{2}(f(1,1)) \geq 3 > 2 \cdot
\text{ord}_{2}(\frac{\partial f}{\partial y}(1, 1)). $ So by Hensel's
lemma (see [Sil1, p.322]), $ f(x, y) $ has a root in $ \Q_{2} \times
\Q_{2}, $ and so $ (\varepsilon (p + q), -pq)_{\Q_{2}} = 1 $ (see
[Weib, Examp.6.2.2, p.253]). Therefore,
$$ \omega _{2} = \sigma _{\varphi }(E / \Q_{2}) =
(-1)^{1 + \text{ord}_{2} \sharp \text{coker}\varphi _{2}}. $$
To calculate the integer $ \sharp \text{coker}\varphi _{2} =
\sharp (E^{\prime }(\Q_{2}) / \varphi _{2}(E(\Q_{2}))), $ we use Lemma 3.8
of [Sc, pp.91, 92]. For this, let
$$ z = - \frac{x}{y}, \ \text{and} \ z^{\prime } =
- \frac{x + \varepsilon (p+q) + pq x^{-1}}{y - pq yx^{-2}} =
- \frac{y}{x^{2} - pq}. $$ From the Chapter IV of [Sil1], one has
$ x = \frac{z}{w(z)} $ and $ y = - \frac{1}{w(z)}, $ where
$ w(z) = z^{3} (1 +  \varepsilon (p+q)z^{2} + \cdots ). $ So
\begin{align*} z^{\prime }&= \frac{w(z)}{z^{2} - pq w(z)^{2}} =
\frac{z^{3} (1 +  \varepsilon (p+q)z^{2} + \cdots )}{z^{2} -
pq z^{6}(1 +  \varepsilon (p+q)z^{2} + \cdots)^{2}} \\
&= z(1 +  \varepsilon (p+q)z^{2} + \cdots ) \cdot (1 + pqz^{4}
(1 +  \varepsilon (p+q)z^{2} + \cdots)^{2} + \cdots ) \\
&= z + (\text{terms of higher degree}),
\end{align*}
i.e., the leading coefficient of $ z^{\prime } $ is $ 1. $ So
$ \mid \varphi _{2}^{\prime }(0) \mid _{2}^{-1} = 1 $ (see [Sc, p.92]),
and so by Lemma 3.8 of [Sc, p.91], we get
$$ \sharp \text{coker}\varphi _{2} =
\frac{\mid \varphi _{2}^{\prime }(0) \mid _{2}^{-1} \cdot \sharp
E(\Q_{2})[\varphi _{2}] \cdot c_{2}(E^{\prime })}{c_{2}(E)} =
\frac{\sharp E(\Q_{2})[\varphi _{2}]
\cdot c_{2}(E^{\prime })}{c_{2}(E)}, $$
where $ c_{2}(E) $ and $ c_{2}(E^{\prime }) $ are the Tamagawa
numbers of $ E $ and $ E^{\prime } $ at $ 2, $ respectively, and
$ E(\Q_{2})[\varphi _{2}] = \text{ker} \varphi _{2} = \{O, (0, 0 )\}. $
So by Lemma 2.1 and Lemma 5.2 above, we get
$ \sharp \text{coker}\varphi _{2} = 2 $ or $ 4, $ that is, \\
 $ \sharp \text{coker}\varphi _{2} = 2 $ if one of the following three
 hypotheses holds: \\
(a) \ $ p \equiv 3 (\text{mod} 8); \ $ (b) \ $ \varepsilon = 1 $
and $ p \equiv 1 (\text{mod} 8); \ $ (c) \ $ \varepsilon = -1 $
and $ p \equiv 5 (\text{mod} 8). $ \\
$ \sharp \text{coker}\varphi _{2} = 4 $ if one of the following three
hypotheses holds: \\
(a$^{\prime }$) \ $ p \equiv 7 (\text{mod} 8); \ $ (b$^{\prime }$) \
$ \varepsilon = 1 $ and $ p \equiv 5 (\text{mod} 8); \ $ (c$^{\prime }$) \
$ \varepsilon = -1 $ and $ p \equiv 1 (\text{mod} 8). $ \\
From this the value of $ \sigma _{\varphi }(E / \Q_{2}) $ and hence
$ \omega _{2} $ is obtained. The proof is completed. \quad $ \Box $
\par \vskip 0.2 cm

On the parity conjecture of some special $ E / \Q $ in (1.1) above, we have
\par \vskip 0.2 cm

{\bf Corollary 5.4.} \ Let $ E / \Q $ be the elliptic curve in
(1.1) above. If one of the
following three hypotheses holds: \\
(1) \ $ \varepsilon = 1 $ and
$ p \equiv 5 \ (\text{mod} \ 8); $ \\
(2) \ $ \varepsilon = -1 $ and
$ p \equiv 3, 5 \ (\text{mod} \ 8); $ \\
(3) \ $ \varepsilon = 1, \ p \equiv 3 \ (\text{mod} \ 8) $
and $ q = a_{1}^{2} + a_{2}^{2} $ with $ (a_{1} + \varepsilon _{1})^{2}
+  (a_{2} + \varepsilon _{2})^{2} = a_{3}^{2} $ for some rational
integers $ a_{1}, a_{2}, a_{3} \in \Z $ and some
$ \varepsilon _{1}, \varepsilon _{2} \in \{1, -1\}. $ \\
Then the parity conjecture is true for $ E / \Q, $
i.e., $ \omega _{E} = (-1)^{\text{rank}E(\Q)}. $
\par \vskip 0.1 cm

{\bf Proof.} \ For the cases (1) and (2), by Theorems 1 and 2 of
[QZ1], we have $ \text{rank}E(\Q) = 0, $ and for the case (3),
by Theorem 3 of [QZ1], we have $ \text{rank}E(\Q) = 1. $ Then the
conclusion follows from Theorem 5.3 above. \quad $ \Box $
\par \vskip 0.2 cm

{\bf Remark.} \ As pointed out by an anonymous referee,
the result of these special $ E / \Q $ in Cor.5.4 above also follows by
Monsky's theorem on the 2-parity conjecture, because their
$ \amalg\hskip-7pt\amalg(E/\Q)[2] $
have been shown to be trivial in [QZ1, Theorems 1,2].
\par \vskip 0.2 cm

{\bf Theorem 5.5.} \ Let $ E / \Q $ be the elliptic curve in
(1.1) and let $ K = \Q(\sqrt{\mu D}) $ be the quadratic number
field with $ D $ in (1.2) and $ \mu = \pm 1. $ We assume that
$ D \equiv \mu \ (\text{mod} \ 4). $ Let $ L(E / \Q, s) =
\Sigma _{n = 1}^{\infty } a_{1}(n) n^{-s} $ be the $ L-$function
as above. Let $ E_{\mu D} / \Q $  be the quadratic $ (\mu D)-$twist
of $ E / \Q, $ and $ \chi _{K} $ be the quadratic Dirichlet character
associated to $ K. $ \\
(1) \ Assume one of the following two hypotheses holds: \\
(a) \ $ \varepsilon = 1 $ and $ p \equiv 5, 7 \ (\text{mod} \ 8); $ \\
(b) \ $ \varepsilon = -1 $ and $ p \equiv 3, 5 \ (\text{mod} \ 8). $ \\
Then $ L(E / \Q, 1) = 2 \Sigma _{n = 1}^{\infty } \frac{a_{1}(n)}{n}
e^{-n\pi /2\sqrt{2pq}}. $ \\
further, for all integer $ r \geq 0, $
\begin{align*}
&L^{(r)}(E / \Q, 1) = 2 \pi \Sigma _{n = 1}^{\infty } a_{1}(n)
\int _{1 /4 \sqrt{2pq}}^{\infty } [\log ^{r} t +
(-1)^{r}\log ^{r}(2^{5}pq t)] e^{-2n\pi t} dt. \ \text{also}, \\
&L(E_{\mu D} / \Q, 1) = (1 + \chi _{K}(-2pq)) \cdot
\Sigma _{n = 1}^{\infty } \frac{a_{1}(n)}{n} \chi _{K}(n) \cdot
e^{-n\pi /2D\sqrt{2pq}},
\end{align*}
In particular, if $ \chi _{K}(-2pq) = -1, $ then
$ L(E_{\mu D} / \Q, 1) = 0. $ \\
(2) \ Assume one of the following two hypotheses holds: \\
(a$^{\prime}$) \ $ \varepsilon = 1 $ and $ p \equiv 1, 3 \ (\text{mod} \ 8); $ \\
(b$^{\prime}$) \ $ \varepsilon = -1 $ and $ p \equiv 1, 7 \ (\text{mod} \ 8). $ \\
Then $ L(E / \Q, 1) = 0, $ \\
further, for all integer $ r \geq 0, $
\begin{align*}
&L^{(r)}(E / \Q, 1) = 2 \pi \Sigma _{n = 1}^{\infty } a_{1}(n)
\int _{1 /4 \sqrt{2pq}}^{\infty } [\log ^{r} t +
(-1)^{r + 1}\log ^{r}(2^{5}pq t)] e^{-2n\pi t} dt. \ \text{also}, \\
&L(E_{\mu D} / \Q, 1) = (1 - \chi _{K}(-2pq)) \cdot
\Sigma _{n = 1}^{\infty }  \frac{a_{1}(n)}{n} \chi _{K}(n) \cdot
e^{-n\pi /2D\sqrt{2pq}}.
\end{align*}
In particular, if $ \chi _{K}(-2pq) = 1, $ then $ L(E_{\mu D} / \Q, 1) = 0. $
\par \vskip 0.1 cm

{\bf Proof.} \ Since $ E/ \Q $ is modular (see [TW],[Wi],[BCDT]), the
function $ f_{E}(z) = \sum _{n = 1}^{\infty } a_{1}(n) e^{2\pi inz} $
satisfies the Hecke equation $ f_{E}(z) = -  \omega _{E} N^{-1}z^{-2}
f(- \frac{1}{Nz}), $ and the differential $ f_{E}(z) dz $ is invariant
under the usual modular group $ \Gamma _{0}(N), $ where $ N = 2^{5}pq $
is the conductor, and $ \omega _{E} $ is the root number of $ E/ \Q. $
Also by assumption, the discriminant $ d(K) = \mu D $ satisfying
$ (d(K), 2 N_{E}) = 1. $ So $ L(E_{\mu D} / \Q, 1) = L(E / \Q, \chi _{K}, 1). $
Hence by Theorem 9.3 of [M, P.61], we have
\begin{align*}  &L(E / \Q, 1) = (1 + \omega _{E}) \Sigma _{n = 1}^{\infty }
\frac{a_{1}(n)}{n} e^{-2n\pi /\sqrt{N}}, \\
&L^{(r)}(E / \Q, 1) = 2 \pi \Sigma _{n = 1}^{\infty } a_{1}(n)
\int _{1 /\sqrt{N}}^{\infty } [\log ^{r} t + \omega _{E}
(-1)^{r}\log ^{r}(N t)] e^{-2n\pi t} dt, \\
&L(E_{\mu D} / \Q, 1) = \Sigma _{n = 1}^{\infty } \frac{a_{1}(n)}{n}
[\chi _{K}(n) + \overline{\chi _{K}}(n) \cdot
\frac{g(\chi _{K})}{g(\overline{\chi _{K}})}\cdot \chi _{K}(-n) \cdot
\omega _{E}] e^{- 2n\pi /\sqrt{N} d(K)},
\end{align*}
where $ g(\chi _{K}) = \sum _{b \ \text{mod} \ d(K)} \chi _{K}(b)
e^{2\pi i b/ d(K)} $ is the Gaussian sum.
Note that $ \chi _{K} (n) = 0, \pm 1 \ (\forall n \in \Z), $ so
$ \overline{\chi _{K}} = \chi _{K}, $ and
$ g(\chi _{K}) = g(\overline{\chi _{K}}). $ Then by our results about
the root numbers in Lemma 5.1 and Theorem 5.3 above, the conclusion
follows.  \quad $ \Box $
\par \vskip 0.2 cm

{\bf Example 5.6.} \ For the elliptic curves
$ E : y^{2} = x (x + 3\varepsilon ) (x + 5\varepsilon ) $ and the
field $ K = \Q (\sqrt{-119}), $ the conductor $ N_{E} =
2^{5} \cdot 3 \cdot 5 = 480 $ and the discriminant $ d(K) = -119. $
By Theorem 5.3 above, the root number of $ E / \Q $ is
$ \omega _{E} = - \varepsilon. $ So for the $ L-$function
$ L(E /\Q, s), $ we have $ L(E /\Q, 1) = 0 $ in the case
$ \varepsilon = 1. $ And in this case, the Mordell-Weil group
$ E(\Q) \cong \Z \times \Z / 2 \Z \times \Z / 2 \Z. $
For the other case $ \varepsilon = -1, \ E(\Q) \cong
\Z / 2 \Z \times \Z / 2 \Z $ (see [QZ1, p.1373]), and by Theorem.5.5
above, $ L(E / \Q, 1) = 2 \Sigma _{n = 1}^{\infty } \frac{a_{1}(n)}{n}
e^{-n\pi /2\sqrt{30}}. $ Moreover, $ d(K) = -119 \equiv 61^{2} \
(\text{mod} \ 4 N_{E}). $ So the
Heegner hypothesis holds for $ E $ and $ K, $ and then there is
a Heegner point $ P_{K} \in E(K) $ such that $ \sigma (2 P_{K}) =
- 2 \omega _{E} P_{K} $ (see [Kol3,4]) because
$ E(\Q)_{\text{tors}} \cong \Z / 2\Z \times \Z / 2\Z $,
where $ \sigma $ is the generator
of the Galois group $ \text{Gal}(K /\Q). $ Since
$ \omega _{E} = - \varepsilon , $ we have $ \sigma (2 P_{K}) =
2 \varepsilon P_{K}. $ Now for any prime number $ l >  37, $
the Galois representation $ \rho _{l} $ is irreducible (see [Cha, p.175]).
Also every such prime number $ l$ satisfies
$ l \nmid d(K), l^{2} \nmid N_{E}, $ so by Cha's theorem in [Cha],
we have $ \text{ord}_{l} \sharp \amalg\hskip-7pt\amalg(E/K) \leq 2
\cdot \text{ord}_{l} ([E(K) : \Z P_{K}]). $ \quad $ \Box $
\par  \vskip 0.2cm

{\bf Remark.} \ This paper is a revised version of the early one ([Q2], 2015).
I thank the anonymous expert for pointing out that the result of Corollary 5.4 above also
follows by Monsky's theorem on the 2-parity conjecture. Some further application toward verifying 
the BSD for a family of elliptic curves will be discussed in a separate paper.
\par  \vskip 0.3cm

{ \bf Acknowledgments }
I would like to thank the referee for helpful suggestions and comments.

\par  \vskip 0.3 cm

\hspace{-0.8cm} {\bf References }
\begin{description}

\item[[BCDT]] C.Breuil, B.Conrad, F.Diamond, R.Taylor, On the
modularity of elliptic curves over $ \Q: $ wild $ 3-$adic exercises,
J.Amer.Math.Soc., 14 (2001), 843-939.

\item[[BS]] M.Bahargava, C.Skinner, A positive proportion of elliptic
curves over $ \Q $ have rank one, arXiv: 1401.0233, 2014.

\item[[BSZ]] M.Bahargava, C.Skinner, W.Zhang, A majority of elliptic
curves over $ \Q $ satisfy the Birch and Swinnerton-Dyer conjecture,
arXiv: 1407.1826 v2, 2014.

\item[[Cha]] B.Cha, Vanishing of some cohomology groups and bounds
for the Shafarevich-Tate groups of elliptic curves, J. Number Theory,
111(2005), 154-178.

\item[[D]] T.Dokchitser, Notes on the parity conjecture,
arXiv: 1009.5389 v2, 2012.

\item[[DD]] T.Dokchitser, V.Dokchitser, Root numbers and parity of ranks
of elliptic curves, arXiv: 0906.1815 v1, 2009.

\item[[DFK]] C.David, J.Fearnley, H.Kisilevsky, On the vanishing of
twisted $ L-$functions of elliptic curves, Experimental Math.,
13 (2004), 185-198.

\item[[Gr]] R.Greenberg, Iwasawa theory for elliptic curves.
In Arithmetic Theory of Elliptic Curves, Lecture Notes in
Math., Vol.1716. New York: Springer-Verlag, 1999, 51-144.

\item[[GV]] R.Greenberg, V.Vatsal, On the Iwasawa invariants of
elliptic curves, Invent. math., 142 (2000), 17-63.

\item[[Kn]] A.W.Knapp, Elliptic Curves, Mathematical Notes 40,
Princeton: Princeton University Press, 1992.

\item[[Kol1]] V.A.Kolyvagin, Finiteness of $ E(\Q) $ and
$ \amalg\hskip-6pt\amalg(E/\Q)$ for a subclass of Weil curves,
(Russian) Izv. Akad. Nauk SSSR Ser. Mat. 52 (1988), 522-540,
670-671; translation in Math. USSR-Izv. 32 (1989), 523-541.

\item[[Kol2]] V.A.Kolyvagin, The Mordell-Weil and Shafarevich-Tate
groups for Weil elliptic curves, (Russian) Izv. Akad. Nauk SSSR Ser.
Mat. 52 (1988), 1154-1180, 1327; translation in Math. USSR-Izv. 33
(1989), 473-499.

\item[[Kol3]] V.A.Kolyvagin, On the Mordell-Weil group and Shafarevich-Tate
group of modular elliptic curves, Proceedings of the international congress of
mathematicians, Kyoto, Japan (1990), 429-436.

\item[[Kr]] K.Kramer, Arithmetic of elliptic curves upon quadratic
extension, Transactions of the American Mathematical Society, 264
(1981), 121-135.

\item[[KT]] K.Kramer, J.Tunnell, Elliptic curves and local
$ \varepsilon-$factors, Compositio Math., 46 (1982), 307-352.

\item[[Ku]] M.Kurihara, On the Tate Shafarevich groups over
cyclotomic fields of an elliptic curve with supersingular
reduction I, Invent. math., 149 (2002), 195-224.

\item[[L]] S.Lang, Algebraic Number Theory, 2nd Edition, New York:
Springer-Verlag, 1994.

\item[[LQ]] F.Li, D.R.Qiu, On several families of elliptic curves with
arbitrary large Selmer groups, Science in China (series A),
53 (2010), 2329-2340.

\item[[M]] J.I.Manin, Cyclotomic fields and modular curves,
Russian Math. Surveys, 26 (1971), 7-78.

\item[[Ma1]] B.Mazur, Rational isogenies of prime degree
(with an appendix by D.Goldfeld), Invent. math., 44 (1978),
129-162.

\item[[Ma2]] B.Mazur, Rational points of Abelian varieties with
values in towers of number fields, Invent. math., 18 (1972),
183-266.

\item[[MR]] B.Mazur, K.Rubin, Ranks of twists of elliptic curves
and Hilbert's tenth problem, Invent. math., 181 (2010),
541-575.

\item[[MSD]] B.Mazur, H.Swinnerton-Dyer, Arithmetic of Weil curves,
Invent. math., 25 (1974), 1-61.

\item[[PR]] B.Perrin-Riou, Arithm\'{e}tique des courbes elliptiques \`{a}
r\'{e}duction supersinguli\`{e}re en $ p, $ Experimental Math.,
12 (2003), 155-186.

\item[[Q1]] D.R.Qiu, On quadratic twists of elliptic curves and
some applications of a refined version of Yu's formula,
Communications in Algebra, 42 (2014), 5050-5064.

\item[[Q2]] D.R.Qiu, On elliptic curves $ y^{2} = x(x + \varepsilon p )(x +
  \varepsilon q) $ and their twists, arXiv: 1511.07581 v1, 2015.

\item[[QZ1]] D.R.Qiu, X.K.Zhang, Mordell-Weil groups and Selmer
groups of twin-prime elliptic curves, Science in China (series A),
45 (2002), 1372-1380.

\item[[QZ2]] D.R.Qiu, X.K.Zhang, Elliptic curves and their torsion
subgroups over number fields of type $ (2, \cdots, 2), $
Science in China (series A), 44 (2001), 159-167.

\item[[Roh]] D.E.Rohrlich, Variation of the root number in families
of elliptic curves, Compositio math., 87 (1993), 119-151.

\item[[Sa]] J.W.Sands, Popescu's conjecture in multi-quadratic extensions,
Contemporary Math., Vol.358, 2004, 127-141.

\item[[Sc]] E.F.Schaefer, Class groups and Selmer groups, J. Number
Theory, 56 (1996), 79-114.

\item[[Se1]] J.-P.Serre, Propri\'{e}t\'{e}s galoisiennes des points d'ordre
fini des courbes elliptiques, Invent. Math., 15 (1972), 259-331.

\item[[Se2]] J.-P. Serre, Local fields, New York: Springer-Verlag,
1979.

\item[[Sil1]] J.H.Silverman, The Arithmetic of Elliptic Curves, GTM
106, New York: Springer-Verlag, 1986.

\item[[Sil2]] J.H.Silverman, Advanced topics in the Arithmetic
of Elliptic Curves, GTM 151, New York: Springer-Verlag, 1999.

\item[[Sk]] C.Skinner, Multiplicative reduction and the cyclotomic
main conjecture for Gl$_{2}, $ arXiv: 1407.1093, 2014.

\item[[SZ]] C.Skinner, W.Zhang, Indivisibility of Heegner points in the
Multiplicative case. arXiv: 1407.1099, 2014.

\item[[St]] M.Stoll, Descent on elliptic curves, arXiv: 0611694 v1,
2006.

\item[[Sz]] K.Szymiczek, $ 2-$ranks of class groups of Witt equivalent
number fields, Annales Mathematicae Silesianae 12 (1998), 53-64.
1979.

\item[[Ta]] J.Tate, Algorithm for determing the type of a singular
fiber in an elliptic pencil, in: Modular functions of one variable,
IV, (Proc. Internat. Summer School, Univ. Antwerp 1972), pp.33-52.
Lecture Notes in Math. 476, Springer, Berlin, 1975.

\item[[TW]] R. Taylor, A.Wiles, Ring-theoretic properties of certain Hecke
algebras, Ann. Math., 141 (1995), 553-572.

\item[[Wa]] L.C.Washington, Introduction to Cyclotomic Fields, 2nd Edition,
New York: Springer-Verlag, 1997.

\item[[Weib]] C.A.Weibel, The K-Book, An Introduction to Algebraic
K-Theory, AMS, Providence, Rhode Island, 2013.

\item[[Wi]] A.Wiles, Modular elliptic curves and Fermat's last theorem,
Ann. Math., 141 (1995), 443-551.

\item[[Zh]] W.Zhang, Selmer groups and the Indivisibility of
Heegner points, Cambridge J. Math., 2 (2014), 191-253.

\end{description}

\end{document}